# The dynamics around stable and unstable Hamiltonian relative equilibria

Juan-Pablo Ortega[1] and Tudor S. Ratiu[2]

**Abstract**

For a symmetric Hamiltonian system, lower bounds for the number of relative equilibria surrounding stable and formally unstable relative equilibria on nearby energy levels are given.

## 1  Introduction

The search for relative equilibria in the presence of nondegeneracy hypotheses has been an extremely active field of research [MO97, RdSD97, MR97, LS98, OR97, O98, CLOR02, He01] during the last few years. In this paper, we will study in a differentiated manner the existence of relative equilibria around *stable* and *formally unstable* equilibria and relative equilibria. We will give estimates on the number of these solutions in terms of readily computable quantities, in order to facilitate the application of these results to specific systems.

A major difference between the bifurcation and persistence results presented in this paper and those in [MO97, RdSD97, MR97, LS98, OR97, O98, He01] is that in our case the solutions obtained are parametrized by energy and not by momentum and, most importantly, our hypotheses do not require the non degeneracy conditions present in all those papers. Consequently our results, particularly theorems 4.1 and 7.3, can be seen as statements on not mere persistence of dynamical elements but on genuine bifurcation phenomena.

The contents of the paper and, in particular the main results, are structured as follows:

- Section 2 contains some preliminaries on symmetric Hamiltonian systems and critical point theory that will be needed in the statements and proofs of the main results.

- Section 3 contains a result (Theorem 3.1) which provides a lower bound for the number of relative equilibria surrounding a stable symmetric Hamiltonian equilibrium whenever a velocity satisfying certain hypotheses can be found.

- Section 4: the superposition of the methods used in Theorem 3.1 with the standard Lyapunov–Schmidt reduction procedure, as well as other techniques dealing with the bifurcation theory of gradient systems, provide in Theorem 4.1 an existence result on branches relative equilibria surrounding formally unstable equilibria.

- Section 5 contains two examples that illustrate the implementation of Theorem 4.1.

- Section 6 is a brief exposition of the Marle–Guillemin–Sternberg normal form [Mar85], [GS84] and the reconstruction equations [O98, RWL99] needed in the next section. The expert can skip this section.

[1]Institut Nonlinéaire de Nice, UMR 129, CNRS-UNSA, 1361, route des Lucioles, 06560 Valbonne, France. Juan-Pablo.Ortega@inln.cnrs.fr.

[2]Institut de Mathématiques Bernoulli, École Polytechnique Fédérale de Lausanne, CH-1015 Lausanne, Switzerland. Tudor.Ratiu@epfl.ch.

- Section 7 presents as main results theorems 7.1 and 7.3, which are the natural generalizations of theorems 3.1 and 4.1, respectively, to the study of relative equilibria surrounding a genuine *relative* equilibrium, using the normal form theory and the reconstruction equations presented in the previous section.

## 2 Preliminaries

$G$–**Hamiltonian systems.** In this paper we will work in the category of symmetric Hamiltonian systems (see, for instance, [AM78]). This means that one considers triples $(M, \omega, h)$, where $\omega$ is a symplectic two–form on the manifold $M$ and $h \in C^\infty(M)$ is a smooth function, called the **Hamiltonian**. Then one associates to $h$ a **Hamiltonian vector field** $X_h$ via the **Hamilton equations**

$$\mathbf{i}_{X_h} \omega = \mathbf{d}h.$$

The symmetries of the system are defined by the left action of a Lie group $G$ on the manifold $M$ that preserves both the symplectic structure $\omega$, that is, the group action is **canonical**, and the Hamiltonian function $h$. The action of $g \in G$ on $m \in M$ will be usually denoted by $g \cdot m$, the space of $G$–invariant smooth functions on $M$ is denoted by $C^\infty(M)^G$, $\mathfrak{g}$ is the Lie algebra of $G$, $\mathfrak{g}^*$ is its dual, and $\exp : \mathfrak{g} \to G$ denotes the exponential map. In most cases we will assume that the $G$–action is also proper and **globally Hamiltonian**, that is, we can associate to it an equivariant **momentum map** $\mathbf{J} : M \to \mathfrak{g}^*$ defined by

$$\mathbf{i}_{\xi_M} \omega = \mathbf{d}\mathbf{J}^\xi,$$

where $\xi_M(m) := (d/dt) \exp t\xi \cdot m|_{t=0}$ is the infinitesimal generator vector field associated to $\xi \in \mathfrak{g}$ and $\mathbf{J}^\xi = \langle \mathbf{J}, \xi \rangle$ is the $\xi$–component of the momentum map $\mathbf{J}$. By **Noether's Theorem**, $\mathbf{J}$ is preserved by the flow of any Hamiltonian vector field associated to any $G$–invariant Hamiltonian function $h \in C^\infty(M)^G$. In particular, the level sets of $\mathbf{J}$ are invariant by the flow of $X_h$.

In the first sections of the paper we will work on a Hamiltonian symplectic vector space $(V, \omega)$, where there is a compact Lie group $G$ acting linearly and canonically. Any such action has an associated equivariant momentum map $\mathbf{J} : V \to \mathfrak{g}^*$ defined by

$$\langle \mathbf{J}(v), \eta \rangle = \frac{1}{2} \omega(\eta \cdot v, v), \qquad \text{for any} \qquad v \in V, \eta \in \mathfrak{g}.$$

The symbol $\eta \cdot v$ denotes the representation of $\mathfrak{g}$ on $V$, which equals $\eta_V(v)$, the value at $v$ of the infinitesimal generator $\eta_V$.

A **relative equilibrium** of the $G$–invariant Hamiltonian $h$ is a point $m \in M$ such that the integral curve $m(t)$ of the Hamiltonian vector field $X_h$ starting at $m$ equals $\exp(t\xi) \cdot m$ for some $\xi \in \mathfrak{g}$. Any such $\xi$ is called a **velocity** or **generator** of the relative equilibrium $m$. Note that if $m$ has a non trivial isotropy subgroup $G_m$, $\xi$ is not uniquely determined. Note also that the $G$–equivariance of the flow of $X_h$ implies that if $m$ is a relative equilibrium with velocity $\xi$ then $g \cdot m$ is also a relative equilibrium but with velocity $\mathrm{Ad}_g \xi$ for any $g \in G$, where $\mathrm{Ad}_g$ is the adjoint representation of $G$ on $\mathfrak{g}$. Thus, we are led to introduce the notion of **distinct relative equilibria**: we say that two relative equilibria are **distinct** when the associated equilibria in the quotient space $M/G$ are distinct. More generally, if $H$ is a closed subgroup of $G$, we say that two relative equilibria are $H$–**distinct** when the associated equilibria in the quotient space $M/H$ are distinct. The topological space $M/G$ is not a manifold in general and the equilibrium needs to be understood in terms of the induced flow on the quotient, that is, an equilibrium in $M/G$ is a point $[m] \in M/G$ such that the quotient flow leaves it fixed.

A key property of symmetric Hamiltonian systems that will be heavily used in this paper is the fact that a point $m \in M$ is a relative equilibrium with velocity $\xi$ if and only if it is a critical point of the so called **augmented Hamiltonian** $h^\xi := h - \mathbf{J}^\xi$. Thus $m \in M$ is a relative equilibrium of the Hamiltonian system with symmetry $(M, \omega, h, G, \mathbf{J})$ with velocity $\xi \in \mathfrak{g}$ if and only if $\mathbf{d}h^\xi(m) = 0$.

If $f \in C^\infty(M)^G$ has a critical point $m$ then $g \cdot m$ is also a critical point of $f$ for any $g \in G$. We shall call **critical orbits** of $f$ the $G$–orbits all of whose points are critical points of $f$.

**Lusternik–Schnirelman category.** For future reference we state the following results:

**Proposition 2.1** *Let $M$ be a compact $G$–manifold, with $G$ a compact Lie group. Any $G$–invariant smooth function $f \in C^\infty(M)^G$ has at least $\mathrm{Cat}(M/G)$ critical orbits.*

In the previous statement, the symbol Cat denotes the Lusternik–Schnirelman category of the quotient compact topological space $M/G$ (the action of $G$ on $M$ does not need to be free and, consequently, the quotient $M/G$ is not in general a manifold). Recall that the Lusternik–Schnirelman category of a compact topological space $M$ is the minimal number of closed contractible sets needed to cover $M$. This result is due to A. Weinstein [W77]. Even though in that reference the result is stated for the case $G = S^1$, the proof that the author provides goes through without any effort in the case in which $G$ is an arbitrary compact Lie group.

Another approach to the search of critical orbits of symmetric functions is the use of the $G$–Lusternik–Schnirelman category, introduced in different versions and degrees of generality by Fadell [Fa85], Clapp and Puppe [CP86, CP91], and Marzantowicz [Mar89]. The equivariant Lusternik–Schnirelman category is not the standard category of the orbit space, which is what we used in the previous paragraphs, but the minimal cardinality of a covering of the $G$–manifold $M$ by $G$–invariant closed subsets that can be equivariantly deformed to an orbit. The use of this definition has allowed Bartsch [Ba94] to provide the following estimate:

**Proposition 2.2** *Let $G$ be a compact Lie group that contains a maximal torus $T$ and that acts linearly on the vector space $V$. Suppose that the vector subspace $V^T$ of $T$–fixed vectors on $V$ is trivial, that is, $V^T = \{0\}$, then any $G$–invariant function defined in the unit sphere of $V$ (the unit sphere of $V$ is defined with the aid of a $G$–invariant norm on it) has at least*

$$\frac{\dim V}{2(1 + \dim G - \dim T)} = \frac{\dim V}{2(1 + \dim G - \mathrm{rank}\, G)}$$

*critical orbits.*

We recall for future use that a $G$–invariant function $f \in C^\infty(M)^G$ on a $G$–space $M$ is called of **Morse–Bott type** when all its critical points $z \in M$ satisfy that $\ker \mathbf{d}^2 f(z) = \mathfrak{g} \cdot z$.

**The Splitting Lemma.** The proof of the following standard result can be found, for instance, in [BrL75].

**Lemma 2.3** *Let $f \in C^\infty(V \times W)$ with $V$ and $W$ finite dimensional vector spaces and such that the mapping $f|_W$, defined by $f|_W(w) := f(0, w)$, has a non–degenerate critical point at $0$. Then there is a local diffeomorphism defined around the point $(0, 0)$, of the form $\psi(v, w) = (v, \psi_1(v, w))$, such that*

$$(f \circ \psi)(v, w) = \bar{f}(v) + Q(w),$$

*where $Q$ is the non–degenerate quadratic form $Q = \frac{1}{2}\mathbf{d}^2 f|_W(0)$, and $\bar{f}$ is a smooth function on $V$.*

## 3   Relative equilibria around a stable equilibrium

In this section we will prove the existence, under certain hypotheses, of relative equilibria around a symmetric stable equilibrium of the system $(V, \omega, h, G, \mathbf{J})$, where $G$ is a compact Lie group that acts canonically and linearly on the symplectic vector space $V$. As we will see in Section 6 (see Remark 6.2) working in the category of linear symplectic spaces implies no loss of generality.

**Theorem 3.1** *Let $(V, \omega, h, G, \mathbf{J})$ be a Hamiltonian $G$–vector space, with $G$ a compact Lie group. Suppose that $h(0) = 0$, $\mathbf{d}h(0) = 0$, and the quadratic form $Q := \mathbf{d}^2 h(0)$ on $V$ is definite. Let $\xi \in \mathfrak{g}$ be such that the quadratic form $\mathbf{d}^2 \mathbf{J}^\xi(0)$ is non degenerate. Then, for each energy value $\epsilon$ small enough, there are at least*

$$\mathrm{Cat}\left(h^{-1}(\epsilon)/G^\xi\right) = \mathrm{Cat}\left(Q^{-1}(\epsilon)/G^\xi\right) \tag{3.1}$$

$G^\xi$–*distinct relative equilibria in* $h^{-1}(\epsilon)$ *whose velocities are (real) multiples of* $\xi$. *The symbol* $G^\xi := \{g \in G \mid \mathrm{Ad}_g \xi = \xi\}$ *denotes the adjoint isotropy of the element* $\xi \in \mathfrak{g}$ *and* Cat *is the Lusternik–Schnirelman category. If the compact Lie group* $G^\xi$ *has a maximal torus* $T^\xi$ *such that the set* $V^{T^\xi}$ *of* $T^\xi$*–fixed vectors on* $V$ *is trivial, that is,* $V^{T^\xi} = \{0\}$, *then there are at least*

$$\frac{\dim V}{2(1 + \dim G^\xi - \mathrm{rank}\, G^\xi)}. \tag{3.2}$$

*distinct relative equilibria in* $h^{-1}(\epsilon)$ *whose velocities are (real) multiples of* $\xi$.

**Remark 3.2** The estimate (3.1) guarantees the existence of at least one relative equilibrium on each nearby level set of the Hamiltonian, since the Lusternik–Schnirelman category of a compact topological space is always at least one. ♦

**Remark 3.3** The hypotheses on the Hamiltonian function, namely $\mathbf{d}h(0) = 0$ and the definiteness of the quadratic form $\mathbf{d}^2 h(0)$, guarantee that the origin is a stable equilibrium of the Hamiltonian vector field $X_h$ (see, for instance, [AM78]). ♦

**Remark 3.4** The optimal way to apply the theorem consists of studying the estimate that it provides in the fixed point spaces of the various isotropy subgroups of the symmetries in the problem. To be more specific, let $(V, \omega, h, G, \mathbf{J})$ be a Hamiltonian system with symmetry with $G$ a compact Lie group. Let $H \subset G$ be an isotropy subgroup of the $G$–action on $V$. It can be easily shown that the vector subspace $V^H$ of $H$–fixed vectors is a symplectic subspace of $V$ and that it is left invariant by the flow associated to $G$–invariant Hamiltonians. Moreover, if $N(H)$ is the normalizer of $H$ in $G$, the group $L := N(H)/H$ acts naturally and canonically on $V^H$ and has associated momentum map $\mathbf{J}_L : V^H \to \mathfrak{l}^*$ given by

$$\mathbf{J}_L(v) = \Lambda^*(\mathbf{J}(v)), \tag{3.3}$$

where $v \in V^H$ and $\Lambda^*$ is the natural $L$–equivariant isomorphism

$$\Lambda^* : (\mathfrak{h}^\circ)^H \longrightarrow \mathfrak{l}^*,$$

between the $H$–fixed point set of vectors in the annihilator of $\mathfrak{h}$ in $\mathfrak{g}^*$ and the dual of the Lie algebra of $L = N(H)/H$ (see [O98, OR02] for the details).

If instead of applying the previous results to the system $(V, \omega, h, G, \mathbf{J})$ we do it on the family of systems $(V^H, \omega|_{V^H}, h|_{V^H}, N(H)/H, \mathbf{J}_L)$ parameterized by the isotropy subgroups $H$ we will obtain more solutions of the problem and, at the same time, we will obtain an estimate on their isotropies (this is especially sharp when we focus on the maximal isotropy subgroups of the action). ♦

**Proof.** Since $\mathbf{d}^2 h(0)$ is definite, the Morse Lemma (see, for instance, [Mil69]) implies that for all $\epsilon$ small enough, the level sets $h^{-1}(\epsilon)$ are compact submanifolds diffeomorphic to spheres. Since $h$ is $G$–invariant these level sets $h^{-1}(\epsilon)$ are also $G$–invariant. At the same time notice that the equivariance of the momentum map $\mathbf{J}$ implies that $\mathbf{J}^\xi \in C^\infty(V)$ is $G^\xi$–invariant and therefore, by Proposition 2.1, the restriction of $\mathbf{J}^\xi$ to the level sets $h^{-1}(\epsilon)$ has at least Cat$\left(h^{-1}(\epsilon)/G^\xi\right)$ critical $G^\xi$–orbits. Let $v(\epsilon)$ be one of those critical points. By the Lagrange Multiplier Theorem (see, for instance, [AMR99, page 211]) there exists a real number (a **multiplier**) $\Lambda(v(\epsilon)) \in \mathbb{R}$ such that

$$\mathbf{dJ}^\xi(v(\epsilon)) = \Lambda(v(\epsilon))\mathbf{d}h(v(\epsilon)). \tag{3.4}$$

The non degeneracy of $\mathbf{d}^2 \mathbf{J}^\xi(0)$ implies that zero is an isolated critical point of $\mathbf{J}^\xi$ hence, by taking $\epsilon$ small enough, we can force the set $\{v \in V \mid h(v) \leq \epsilon\}$ (whose boundary is the level set $h^{-1}(\epsilon)$) to contain only zero as a critical point of $\mathbf{J}^\xi$. If we restrict $\epsilon$ to that range, we can guarantee that the multiplier $\Lambda(v(\epsilon))$ in (3.4) is not zero since otherwise $v(\epsilon)$ would be a critical point of $\mathbf{J}^\xi$ in $\{v \in V \mid h(v) \leq \epsilon\}$ which is impossible by construction. This circumstance and the linearity of $\mathbf{J}^\xi$ in $\xi$ implies that we can rewrite (3.4) as

$$\mathbf{d}\left(h - \mathbf{J}^{\xi/\Lambda(v(\epsilon))}\right)(v(\epsilon)) = 0,$$

that is, the point $v(\epsilon)$ is a relative equilibrium of the vector field $X_h$ with velocity $\xi/\Lambda(v(\epsilon))$.

The fact that
$$\operatorname{Cat}\left(h^{-1}(\epsilon)/G^\xi\right) = \operatorname{Cat}\left(Q^{-1}(\epsilon)/G^\xi\right).$$

is a consequence of the equivariant Morse Lemma (see [Bott82] and the Appendix of [VvdM95]) by virtue of which there exists a local $G$–equivariant diffeomorphism $\psi$ of $V$ around the origin such that $h \circ \Psi = Q$. Since the Lusternik–Schnirelman category is a topological invariant, the equality follows. The second estimate (3.2) follows from Proposition 2.2. ∎

## 4 Relative equilibria around formally unstable equilibria

In this section we will present a result concerning the bifurcation of relative equilibria from a formally unstable equilibrium. The motivation for this result comes after realizing that the stability hypothesis in the statement of Theorem 3.1 is too strong. We illustrate this fact by giving a very simple example in which the hypotheses of Theorem 3.1 are violated due to the absence of the definiteness hypothesis and nevertheless there exist relative equilibria around the equilibrium in question. Let $V = \mathbb{R}^4$ endowed with the symplectic structure $\omega = \mathbf{d}q_1 \wedge \mathbf{d}p_1 + \mathbf{d}q_2 \wedge \mathbf{d}p_2$. Consider the canonical action of the group $S^1$ given by
$$(e^{i\theta}, (q_1, q_2, p_1, p_2)) \longmapsto (R_\theta(q_1, p_1), R_\theta(q_2, p_2)),$$

where $R_\theta(q_i, p_i)$ denotes the rotation with angle $\theta$ of the vector $(q_i, p_i)$. This action has an equivariant momentum map $\mathbf{J} : \mathbb{R}^4 \to \mathbb{R}$ associated given by $\mathbf{J}(q_1, q_2, p_1, p_2) = \frac{1}{2}(q_1^2 + p_1^2 - q_2^2 - p_2^2)$. Consider now the $S^1$–invariant Hamiltonian
$$h(q_1, q_2, p_1, p_2) = (q_1^2 + p_1^2) - 2(q_2^2 + p_2^2) + (q_1^2 + p_1^2)(q_2^2 + p_2^2).$$

Clearly the definiteness hypothesis in Theorem 3.1 does not hold for $h$. Nevertheless, since
$$\mathbf{d}(h - \mathbf{J}^\xi)(q_1, q_2, p_1, p_2) = (q_1(2 + 2(p_2^2 + q_2^2) - \xi), p_1(2 + 2(p_2^2 + q_2^2) - \xi),$$
$$q_2(\xi - 4 + 2(p_1^2 + q_1^2)), p_2(\xi - 4 + 2(p_1^2 + q_1^2))),$$

any point of the form $(0, q_2, 0, p_2)$ is a $S^1$–relative equilibrium with velocity $\xi = 4$. The same can be said about the points of the form $(q_1, 0, p_1, 0)$, with velocity $\xi = 2$.

The following result is capable of predicting these critical elements. More explicitly, we will show that even if $\mathbf{d}^2 h(0)$ is indefinite, under certain circumstances, the existence of relative equilibria around a given equilibrium is guaranteed.

**Theorem 4.1** *Let $(V, \omega, h, G, \mathbf{J})$ be a Hamiltonian $G$–vector space, with $G$ a compact Lie group. Suppose that $h(0) = 0$ and $\mathbf{d}h(0) = 0$. Let $\xi \in \mathfrak{g}$ be a root of the polynomial equation:*
$$\det\left(\mathbf{d}^2\left(h - \mathbf{J}^\xi\right)(0)\right) = 0. \tag{4.1}$$

*Define*
$$V_0 := \ker\left(\mathbf{d}^2\left(h - \mathbf{J}^\xi\right)(0)\right)$$

*and suppose that:*

**(i)** *The restriction of the quadratic form $Q := \mathbf{d}^2 h(0)|_{V_0}$ on $V_0$ is definite.*

**(ii)** *Let $\|\cdot\|$ be the norm on $V_0$ defined by $\|v_0\| := \mathbf{d}^2 h(0)(v_0, v_0)$, $v_0 \in V_0$. This map is indeed a norm due to the definiteness assumption on $\mathbf{d}^2 h(0)|_{V_0}$ (if $\mathbf{d}^2 h(0)|_{V_0}$ is negative definite, a minus sign is needed in the definition). Let $l = \dim V_0$ and $S^{l-1}$ be the unit sphere in $V_0$. The function $j \in C^\infty(S^{l-1})$ defined by $j(u) := \frac{1}{2}\mathbf{d}^2\mathbf{J}^\xi(0)(u, u)$, $u \in S^{l-1}$, is of Morse–Bott type with respect to the $G^\xi$–action on $S^{l-1}$.*

*Then, there are at least*

$$\mathrm{Cat}\left(h|_{V_0}^{-1}(\epsilon)/G^\xi\right) = \mathrm{Cat}\left(Q^{-1}(\epsilon)/G^\xi\right) \qquad (4.2)$$

$G^\xi$*–distinct relative equilibria of h on each of its energy levels near zero. The velocities of these relative equilibria are close to* $\xi$*. The symbol* $G^\xi$ *denotes the adjoint isotropy of the element* $\xi \in \mathfrak{g}$ *and* Cat *the Lusternik–Schnirelman category.*

*If the subgroup* $G^\xi$ *contains a maximal torus* $T^\xi$ *such that the subspace* $(V_0)^{T^\xi}$ *of* $T^\xi$*–fixed point vectors in* $V_0$ *is trivial, that is* $(V_0)^{T^\xi} = \{0\}$*, then the estimate (4.2) can be replaced by*

$$\frac{\dim V_0}{2(1 + \dim G^\xi - \dim T^\xi)}. \qquad (4.3)$$

**Remark 4.2** The hypothesis on the function $j$ being Morse–Bott with respect to the $G^\xi$–action on $S^{l-1}$ holds generically [GoMac88].  ♦

Before we proceed to prove the theorem we see how it is actually capable of predicting the relative equilibria that we discussed in the motivational example preceding the statement. Indeed, a straightforward calculation shows that in that case, the equation on $\xi$

$$\det\left(\mathbf{d}^2\left(h - \mathbf{J}^\xi\right)(0)\right) = 0$$

has $\xi = \{2, 4\}$ as roots. We associate to each of these roots the spaces

$$V_0^2 = \{(q_1, 0, p_1, 0) \in V \mid q_1, p_1 \in \mathbb{R}\}, \qquad V_0^4 = \{(0, q_2, 0, p_2) \in V \mid q_2, p_2 \in \mathbb{R}\}.$$

The restriction of $\mathbf{d}^2 h(0)$ to both spaces is definite and the corresponding spheres $Q^{-1}(\epsilon)$ amount to circles on which the symmetry group acts transitively forcing the Morse–Bott hypothesis on the functions $j$ to hold. Consequently, Theorem 4.1 provides us with the relative equilibria that we found by hand in this example.

**Proof.** Let $\mathfrak{g}^{G^\xi}$ be the set of elements in $\mathfrak{g}$ fixed by the adjoint action of the subgroup $G^\xi$ on $\mathfrak{g}$. Note that, by the definition of $G^\xi$, $\xi \in \mathfrak{g}^{G^\xi}$. Let $F : V \times \mathfrak{g}^{G^\xi} \to V$ be the mapping defined by

$$F(v, \alpha) := \nabla_V \left(h - \mathbf{J}^{\xi + \alpha}\right)(v), \qquad v \in V,\ \alpha \in \mathfrak{g}^{G^\xi},$$

where the symbol $\nabla_V$ denotes the gradient defined with the aid of a $G$–invariant inner product on $V$, always available by the compactness of $G$. We will search the relative equilibria of the system by looking for the zeros of the mapping $F$.

***Step 1: Lyapunov–Schmidt Reduction.*** We start this study by first performing a Lyapunov–Schmidt reduction on $F$ (see [GoS85]). Let $L : V \to V$ be the mapping defined by $L(v) = \mathbf{d}F(0, 0) \cdot v$. It is easy to show that for any $v, w \in V$, $\langle L(v), w \rangle = \mathbf{d}^2(h - \mathbf{J}^\xi)(0)(v, w)$ and therefore $V_0 = \ker L$. Notice that due to the $G^\xi$–equivariance of $L$, the subspace $V_0$ is $G^\xi$–invariant. Let $V_1$ be a $G^\xi$–invariant complement to $V_0$ in $V$, that is, $V = V_0 \oplus V_1$. Let $\mathbb{P} : V \to V_0$ be the canonical $G^\xi$–equivariant projection associated to this splitting and $v = v_0 + v_1$ be the decomposition of an arbitrary element $v \in V$ in terms of its $V_0$ and $V_1$ components. The equation $(\mathbb{I} - \mathbb{P})F(v_0 + v_1, \alpha) = 0$ defines, via the Implicit Function Theorem, a $G^\xi$–equivariant mapping $v_1 : V_0 \times \mathfrak{g}^{G^\xi} \to V_1$ such that

$$(\mathbb{I} - \mathbb{P})F(v_0 + v_1(v_0, \alpha), \alpha) = 0. \qquad (4.4)$$

***Step 2: Properties of*** $v_1$**.** The function $v_1$ satisfies the properties that we collect in the following lemma:

**Lemma 4.3** *The function* $v_1$ *defined in (4.4) satisfies the following properties:*

**(i)** $v_1(0, \alpha) = 0$,

**(ii)** $D_{V_0} v_1(0,0) = 0$,

**(iii)** $\mathbf{d}^2(h - \mathbf{J}^\xi)(0)(D_{V_0,\alpha} v_1(0,0) \cdot (v, \alpha), z) = \mathbf{d}^2 \mathbf{J}^\alpha(0)(v, z)$.

for any $\alpha \in \mathfrak{g}^{G_\xi}$, $u, v, w \in V_0$, and $z \in V_1$. The symbols $D_{V_0}$ and $D_{V_0,\alpha}$ denote the partial Fréchet derivatives relative to the $V_0$–variable and the second partial derivative relative to the two variables $V_0$ and $\mathfrak{g}^{G_\xi}$.

**Proof** **(i)** $v_1(0, \alpha) = 0$ satisfies equation (4.4). The uniqueness of the solutions obtained by the Implicit Function Theorem establishes the result.
**(ii)** Equation (4.4) is equivalent to

$$\mathbf{d}(h - \mathbf{J}^{\xi+\alpha})(v_0 + v_1(v_0, \alpha)) \cdot z = 0 \quad \text{for all} \quad v_0 \in V_0,\ z \in V_1,\ \alpha \in \mathfrak{g}^{G_\xi}. \tag{4.5}$$

We take the derivative of this expression relative to $V_0$ at the point $(0,0) \in V_0 \times \mathfrak{g}^{G_\xi}$ in the direction $v \in V_0$. One gets

$$0 = \mathbf{d}^2(h - \mathbf{J}^\xi)(0)(v + D_{V_0} v_1(0,0) \cdot v, z) = \mathbf{d}^2(h - \mathbf{J}^\xi)(0)(D_{V_0} v_1(0,0) \cdot v, z)$$

because $v \in V_0 = \ker L$. Therefore, since $\mathbf{d}^2(h - \mathbf{J}^\xi)(0)$ is non degenerate on $V_1$, it follows that $D_{V_0} v_1(0,0) = 0$.
**(iii)** We take the derivative of expression (4.5) relative to $V_0$ and $\mathfrak{g}^{G_\xi}$ at the point $(0,0) \in V_0 \times \mathfrak{g}^{G_\xi}$ in the directions $v \in V_0$ and $\alpha \in \mathfrak{g}^{G_\xi}$. We obtain

$$\left.\frac{d}{dt}\right|_{t=0} \left.\frac{d}{ds}\right|_{s=0} \mathbf{d}(h - \mathbf{J}^{\xi+s\alpha})(tv + v_1(tv, s\alpha)) \cdot z = 0.$$

Using properties **(i)** and **(ii)** in the computation of this derivative the result follows. ▼

**Step 3: The Bifurcation Equation.** With all these ingredients, the final Lyapunov–Schmidt $G^\xi$–equivariant **reduced equation** is given by $B : V_0 \times \mathfrak{g}^{G_\xi} \to V_0$, where

$$\begin{aligned} B(v_0, \alpha) &= \mathbb{P} F(v_0 + v_1(v_0, \alpha), \alpha) = \mathbb{P} \nabla_V (h - \mathbf{J}^{\xi+\alpha})(v_0 + v_1(v_0, \alpha)) \\ &= \nabla_V (h - \mathbf{J}^{\xi+\alpha})(v_0 + v_1(v_0, \alpha)) \qquad \text{(by (4.4))}. \end{aligned} \tag{4.6}$$

Hence, we have reduced the problem of finding the zeros of $F$ to that of finding the zeros of the $G^\xi$–equivariant map $B$ which is defined in a smaller dimensional space. This reduction technique has already been exploited in the symmetric Hamiltonian framework in [CLOR02, COR02]. As it was also noticed in those references, the reduced equation $B$ is the gradient of a $G^\xi$–invariant function defined on $V_0$, that is,

$$B(v_0, \alpha) = \nabla_{V_0} g(v_0, \alpha),$$

where the function $g : V_0 \times \mathfrak{g}^{G_\xi} \to V_0$ is given by

$$g(v_0, \alpha) = (h - \mathbf{J}^{\xi+\alpha})(v_0 + v_1(v_0, \alpha)).$$

We verify that this is indeed the case. Note first that for any $w \in V_1$ we have that

$$\begin{aligned} \langle F(v_0 + v_1(v_0, \alpha), \alpha), w \rangle &= \langle F(v_0 + v_1(v_0, \alpha), \alpha), (\mathbb{I} - \mathbb{P})w \rangle \\ &= \langle (\mathbb{I} - \mathbb{P}) F(v_0 + v_1(v_0, \alpha), \alpha), w \rangle = 0. \end{aligned} \tag{4.7}$$

Where the last equality follows from the construction of the function $v_1$ through expression (4.4). Now, let $u \in V_0$ arbitrary. We write:

$$\begin{aligned} \langle B(v_0, \alpha), u \rangle &= \langle \mathbb{P} F(v_0 + v_1(v_0, \alpha), \alpha), u \rangle \\ &= \langle F(v_0 + v_1(v_0, \alpha), \alpha), u \rangle \\ &= \langle F(v_0 + v_1(v_0, \alpha), \alpha), u + D_{V_0} v_1(v_0, \alpha) \cdot u \rangle \quad \text{(by (4.7))} \\ &= \langle \nabla_V (h - \mathbf{J}^{\xi+\alpha})(v_0 + v_1(v_0, \alpha)), u + D_{V_0} v_1(v_0, \alpha) \cdot u \rangle \\ &= \mathbf{d}(h - \mathbf{J}^{\xi+\alpha})(v_0 + v_1(v_0, \alpha)) \cdot (u + D_{V_0} v_1(v_0, \alpha) \cdot u) \\ &= \mathbf{d} g(v_0, \alpha) \cdot u = \langle \nabla_{V_0} g(v_0, \alpha), u \rangle, \end{aligned}$$

as required. This construction is a particular case of the one carried out in [CLOR02], [COR02] and [GMSD95].

The following lemma provides two additional properties of the reduced bifurcation equation that will be used later on. The proof is a straightforward differentiation of the function $B$ aided by the properties in Lemma 4.3.

**Lemma 4.4** *The reduced bifurcation equation satisfies the following two properties:*

**(i)** $D_{V_0} B(0,0) = 0$,

**(ii)** $\langle D_{V_0,\alpha} B(0,0)(v_0, \alpha), w_0 \rangle = -\mathbf{d}^2 \mathbf{J}^\alpha(0)(v_0, w_0)$,

*for any $v_0, w_0 \in V_0$ and any $\alpha \in \mathfrak{g}^{G^\xi}$.*

***Step 4: Critical Points and Lagrange Multipliers.*** We now define, for any $\alpha, \beta \in \mathfrak{g}^{G^\xi}$, the functions:
$$H_\alpha(v_0) := h(v_0 + v_1(v_0, \alpha)), \qquad \mathbf{J}^\beta_\alpha(v_0) := \mathbf{J}^\beta(v_0 + v_1(v_0, \alpha)). \tag{4.8}$$

Using the properties in Lemma 4.3 and the fact that $\mathbf{d}h(0) = 0$ it is easy to see that for any $\alpha \in \mathfrak{g}^{G^\xi}$
$$\mathbf{d}H_\alpha(0) = 0 \quad \text{and} \quad \mathbf{d}^2 H_0(0) = \mathbf{d}^2 h(0)|_{V_0}. \tag{4.9}$$

The definiteness hypothesis on $\mathbf{d}^2 h(0)|_{V_0}$ and the invariance properties of $h$ allow us to define a $G$–invariant norm $\|\cdot\|$ on $V_0$ by taking
$$\|v_0\|^2 := \mathbf{d}^2 h(0)(v_0, v_0). \tag{4.10}$$

Moreover, the Splitting Lemma 2.3 and (4.9) guarantee the existence of a local $G^\xi$–equivariant change of variables on $V_0$ around the origin in which the function $H_\alpha$ takes the form
$$H_\alpha(v_0) = \|v_0\|^2 + f(\alpha), \tag{4.11}$$

where $f : \mathfrak{g}^{G^\xi} \to \mathbb{R}$ is a smooth function such that $f(0) = 0$. Note that (4.11) implies that for a fixed value of the parameter $\alpha$, the level sets of the function $H_\alpha$ are $G^\xi$–equivariantly diffeomorphic to spheres provided that we stay close enough to the origin in $V_0$.

We will now follow a strategy similar to the one presented in Theorem 3.1 in order to establish the generic part in the statement of the theorem. For any $\alpha, \beta \in \mathfrak{g}^{G^\xi}$, the mapping $\mathbf{J}^{\xi+\beta}_\alpha \in C^\infty(V_0)$ is $G^\xi$–invariant and therefore, its restriction to the level sets $H_\alpha^{-1}(\epsilon)$, with $\alpha$ and $\epsilon$ as in the previous lemma, has at least
$$\text{Cat}\left(H_\alpha^{-1}(\epsilon)/G^\xi\right) \tag{4.12}$$

critical $G^\xi$–orbits (by Proposition 2.1). Let $v_0(\epsilon, \alpha, \beta)$ be one of those critical points. Again, by the Lagrange Multiplier Theorem [AMR99, page 211], there exists a multiplier $\Lambda(\epsilon, \alpha, \beta) \in \mathbb{R}$ such that
$$\mathbf{dJ}^{\xi+\beta}_\alpha(v_0(\epsilon, \alpha, \beta)) = \Lambda(\epsilon, \alpha, \beta) \mathbf{d}H_\alpha(v_0(\epsilon, \alpha, \beta)). \tag{4.13}$$

***Step 5: The Blow–Up Argument.*** In the following paragraphs we will prove that if we reparametrize the mapping $v_0(\epsilon, \alpha, \beta)$ that describes the "branch" of critical points of $\mathbf{J}^{\xi+\beta}_\alpha$ on the level sets of $H_\alpha$ with the norm of $v_0$ instead of with $\epsilon$, we can choose the resulting function to be smooth. We will denote the norm of $v_0$ by $r$. Recall that by (4.11), the relation between $r$ and $\epsilon$ is given, for a fixed $\alpha$, by $\epsilon = r^2 + f(\alpha)$. Let $v_0(r, \alpha, \beta)$ be the function obtained out of $v_0(\epsilon, \alpha, \beta)$ via that relation. As we just said, we will see that the genericity hypotheses under which we are working will guarantee the local smoothness around the origin of $v_0(r, \alpha, \beta)$. Indeed, let us first reformulate our problem using polar coordinates on $V_0$ (blow-up), that is, $v_0 = ru$, with $r \in \mathbb{R}$ and $u \in S^{l-1}$, $l := \dim V_0$, and $S^{l-1}$ is the unit sphere on $V_0$, defined via the norm (4.10). We now define:
$$\bar{H}_\alpha(r, u) := H_\alpha(ru), \qquad \bar{\mathbf{J}}^{\xi+\beta}_\alpha(r, u) := \mathbf{J}^{\xi+\beta}_\alpha(ru). \tag{4.14}$$

The function $\bar{\mathbf{J}}_\alpha^{\xi+\beta}$ can be rewritten as

$$\bar{\mathbf{J}}_\alpha^{\xi+\beta}(r,u) = r^2 \widehat{\mathbf{J}}_\alpha^{\xi+\beta}(r,u),$$

where

$$\widehat{\mathbf{J}}_\alpha^{\xi+\beta}(r,u) = f_{\alpha,\beta}(u) + g_{\alpha,\beta}(r,u),$$

with $f_{\alpha,\beta}$ and $g_{\alpha,\beta}$ smooth functions on their arguments such that $g_{\alpha,\beta}(0,u) = 0$ for any $u \in S^{l-1}$, $\alpha, \beta \in \mathfrak{g}^{G^\xi}$, and

$$f_{\alpha,\beta}(u) = \frac{1}{2}\mathbf{d}^2\mathbf{J}^{\xi+\beta}(0)(u + D_{V_0}v_1(0,\alpha) \cdot u, u + D_{V_0}v_1(0,\alpha) \cdot u).$$

Since for a fixed value of the parameter $\alpha$, the level sets of $H_\alpha$ are spheres ($r$ is constant), the critical points of $\widehat{\mathbf{J}}_\alpha^{\xi+\beta}|_{H_\alpha^{-1}(\epsilon)}$ coincide with the critical points of $\bar{\mathbf{J}}_\alpha^{\xi+\beta}|_{H_\alpha^{-1}(\epsilon)}$, which is what we are trying to describe.

**Step 6: Smoothness of the branches of critical points.** In order to show that these critical points come in smooth branches, consider the $G^\xi$–invariant function $j$ on the sphere $S^{l-1}$, defined by

$$j(u) := \frac{1}{2}\mathbf{d}^2\mathbf{J}^\xi(0)(u,u), \quad u \in S^{l-1}.$$

Let $u_0 \in S^{l-1}$ be one of its critical orbits provided, for instance, by an estimate of the form (4.12). Due to the $G^\xi$–invariance of $j$, $u_0$ is inevitably a degenerate critical point of $j$. Given that by hypothesis $j$ is a Morse–Bott function with respect to the $G^\xi$–action, we have that,

$$\ker \mathbf{d}^2 j(u_0) = \mathfrak{g}^\xi \cdot u_0,$$

where $\mathfrak{g}^\xi \cdot u_0$ is the tangent space at the point $u_0$ to the $G^\xi$–orbit that goes through it. Let now $\sigma$ be a local cross–section of the homogeneous space $G^\xi/G^\xi_{u_0}$, that is, a differentiable map $\sigma : \mathcal{Z} \to G^\xi$, where $\mathcal{Z}$ is an open neighborhood of $G^\xi_{u_0}$ in the homogeneous space $G^\xi/G^\xi_{u_0}$ such that $\sigma(G^\xi_{u_0}) = e$ and $\sigma(z) \in z$, for $z \in \mathcal{Z}$. The existence of these local cross–sections is well known (see for instance [Che46, page 109]). The Slice Theorem [Pal61, propositions 2.1.2 and 2.1.4] guarantees the existence of a submanifold $\mathcal{S}_{u_0}$ of $S^{l-1}$ going through $u_0$ (the $G^\xi$–slice through $u_0$), such that the product $\mathcal{Z} \times \mathcal{S}_{u_0}$ is diffeomorphic to a neighborhood of $u_0$ in $S^{l-1}$ via the map $(gG^\xi_{u_0}, u) \mapsto \sigma(gG^\xi_{u_0}) \cdot u$. When $\mathfrak{g}^\xi \cdot u_0 = T_{u_0}S^{l-1}$ then $\mathcal{S}_{u_0} = \{u_0\}$ and all subsequent arguments have obvious simplifications. Let $(U, \psi = (\psi_1, \psi_2))$ be a product chart for the product manifold $\mathcal{Z} \times \mathcal{S}_{u_0}$ around the point $(G^\xi_{u_0}, u_0)$ such that $\psi(G^\xi_{u_0}, u_0) = (0,0)$. Denote by $(z,s)$ the elements in $\psi(U)$ that we can use to parametrize a neighborhood of $u_0$ in $S^{l-1}$ via the map $\varphi : \psi(U) \to S^{l-1}$ given by $(z,s) \mapsto \sigma(\psi_1^{-1}(z)) \cdot \psi_2^{-1}(s)$. Notice that $\varphi(0,0) = u_0$ and $\mathfrak{g}^\xi \cdot u_0 = T_{u_0}\varphi(\psi_1(\mathcal{Z}) \times \{0\})$.

We now go back to the description of the critical points of $\widehat{\mathbf{J}}_\alpha^{\xi+\beta}$. Since we are interested on how these critical points behave when we move around $u_0$ we will write the function $\widehat{\mathbf{J}}_\alpha^{\xi+\beta}$ using the diffeomorphism $\varphi$. First, the $G^\xi$–invariance of $\widehat{\mathbf{J}}_\alpha^{\xi+\beta}$ implies that its representative in $(z,s)$ coordinates does not depend on $z$, that is, it has the form

$$\widehat{\mathbf{J}}_\alpha^{\xi+\beta}(r,s) = f_{\alpha,\beta}(s) + g_{\alpha,\beta}(r,s),$$

where $f_{0,0}(s) = j(s)$. Second, since $\mathbf{d}j(u_0) = 0$, then $\mathbf{d}_s\widehat{\mathbf{J}}_0^\xi(0,0) = 0$. Also, since $\mathbf{d}^2 j(u_0)|_{T_{u_0}\mathcal{S}_{u_0}}$ is non degenerate, so is $\mathbf{d}_s^2\widehat{\mathbf{J}}_0^\xi(0,0)$, hence we can define via the Implicit Function Theorem a smooth function $s(r,\alpha,\beta)$ such that the points on $V_0$ of the form $r\varphi(z, s(r,\alpha,\beta))$ constitute critical orbits of the restriction of $\widehat{\mathbf{J}}_\alpha^{\xi+\beta}$ to the level sets of $H_\alpha$, that is,

$$\mathbf{d}_s\hat{\mathbf{J}}_\alpha^{\xi+\beta}(r, s(r,\alpha,\beta)) = 0.$$

Consequently, the smooth branch that we are looking for is:

$$v_0(r,\alpha,\beta) := r\varphi(0, s(r,\alpha,\beta)) = r\psi_2^{-1}(s(r,\alpha,\beta)). \tag{4.15}$$

As a corollary to the preceding ideas we obtain that the Lagrange multiplier $\Lambda(\epsilon, \alpha, \beta) \in \mathbb{R}$ introduced in (4.13) is smooth in its arguments if we reparametrize it as a function of the form $\Lambda(r, \alpha, \beta)$. Indeed, if we pair both sides of (4.13), using the new parameterization, with $v_0(r, \alpha, \beta)$ we have that

$$\Lambda(r, \alpha, \beta) = \frac{\mathbf{dJ}_\alpha^{\xi+\beta}(v_0(r, \alpha, \beta)) \cdot v_0(r, \alpha, \beta)}{\mathbf{d}H_\alpha(v_0(r, \alpha, \beta)) \cdot v_0(r, \alpha, \beta)}.$$

As we can easily deduce by looking at (4.11), the denominator of this expression is different from zero as long as we are not at the origin, that is, when $r = 0$. Elsewhere, the function $\Lambda(r, \alpha, \beta)$ is a combination of smooth objects, thereby smooth. In the following Lemma we see that actually the origin is not a singularity and that the function $\Lambda$ is smooth also in there.

**Lemma 4.5** *Let $\Lambda(r, \alpha, \beta)$ be the multiplier introduced in the previous paragraphs. Then, the function $\Lambda(r, \alpha, \beta)$ is smooth at the point $(0, 0, 0)$ and, moreover we have that:*

$$\Lambda(0, 0, 0) = 1.$$

**Proof.** We will deal with this problem using polar coordinates. Let $\bar{H}_\alpha(r, u)$ and $\bar{\mathbf{J}}_\alpha^{\xi+\beta}(r, u)$ be the functions introduced in (4.14). Recall that

$$\bar{\mathbf{J}}_\alpha^{\xi+\beta}(r, u) = r^2 \left[ \frac{1}{2} \mathbf{d}^2 \mathbf{J}^{\xi+\beta}(0)(u + D_{V_0} v_1(0, \alpha) \cdot u, u + D_{V_0} v_1(0, \alpha) \cdot u) + g_{\alpha, \beta}(r, u) \right],$$

and

$$\bar{H}_\alpha(r, u) = r^2 \left[ \frac{1}{2} \mathbf{d}^2 h(0)(u + D_{V_0} v_1(0, \alpha) \cdot u, u + D_{V_0} v_1(0, \alpha) \cdot u) + q_{\alpha, \beta}(r, u) \right],$$

where $g_{\alpha, \beta}$ and $q_{\alpha, \beta}$ are smooth functions such that $g_{\alpha, \beta}(0, u) = q_{\alpha, \beta}(0, u) = 0$ for any $u \in S^{l-1}$, $\alpha, \beta \in \mathfrak{g}^{G^\xi}$. It is easy to see that

$$\frac{\partial \bar{\mathbf{J}}_\alpha^{\xi+\beta}}{\partial r}(r, u) = 2r \left[ \frac{1}{2} \mathbf{d}^2 \mathbf{J}^{\xi+\beta}(0)(u + D_{V_0} v_1(0, \alpha) \cdot u, u + D_{V_0} v_1(0, \alpha) \cdot u) + g_{\alpha, \beta}(r, u) \right]$$
$$+ r^2 \frac{\partial g_{\alpha, \beta}}{\partial r}(r, u),$$

$$\frac{\partial \bar{H}_\alpha}{\partial r}(r, u) = 2r \left[ \frac{1}{2} \mathbf{d}^2 h(0)(u + D_{V_0} v_1(0, \alpha) \cdot u, u + D_{V_0} v_1(0, \alpha) \cdot u) + q_{\alpha, \beta}(r, u) \right]$$
$$+ r^2 \frac{\partial q_{\alpha, \beta}}{\partial r}(r, u),$$

$$\frac{\partial \bar{\mathbf{J}}_\alpha^{\xi+\beta}}{\partial r}(r, u) = \mathbf{dJ}_\alpha^{\xi+\beta}(ru) \cdot u.$$

We pair the defining expression of the multiplier (4.13) on both sides with $\psi_2^{-1}(s(r, \alpha, \beta))$. By (4.15) and the three relations above we get

$$\Lambda(r, \alpha, \beta) = \frac{\mathbf{dJ}_\alpha^{\xi+\beta}(v_0(r, \alpha, \beta)) \cdot u}{\mathbf{d}H_\alpha(v_0(r, \alpha, \beta)) \cdot u}$$
$$= \frac{2\left[\frac{1}{2}\mathbf{d}^2\mathbf{J}^{\xi+\beta}(0)(u + D_{V_0}v_1(0, \alpha) \cdot u, u + D_{V_0}v_1(0, \alpha) \cdot u) + g_{\alpha,\beta}(r, u)\right] + r\frac{\partial g_{\alpha,\beta}}{\partial r}(r, u)}{2\left[\frac{1}{2}\mathbf{d}^2 h(0)(u + D_{V_0}v_1(0, \alpha) \cdot u, u + D_{V_0}v_1(0, \alpha) \cdot u) + q_{\alpha,\beta}(r, u)\right] + r\frac{\partial q_{\alpha,\beta}}{\partial r}(r, u)}, \quad (4.16)$$

where in the previous expression the symbol $u$ denotes $\psi_2^{-1}(s(r, \alpha, \beta))$ (see 4.15). Notice that since we have had one cancellation of $r$, the previous expression is not singular anymore at the point $(0, 0, 0)$. Moreover,

$$\Lambda(0, 0, 0) = \frac{\mathbf{d}^2\mathbf{J}^\xi(0)(u, u)}{\mathbf{d}^2 h(0)(u, u)} = 1,$$

given that $u \in V_0 = \ker\left(\mathbf{d}^2\left(h - \mathbf{J}^\xi\right)(0)\right)$ and, therefore $\mathbf{d}^2\mathbf{J}^\xi(0)(u,u) = \mathbf{d}^2 h(0)(u,u) \neq 0$, by the definiteness hypothesis on $\mathbf{d}^2 h(0)|_{V_0}$.  ▼

### Step 7: Reduction of the problem to a scalar equation.

**Lemma 4.6** *Let $\Lambda(r, \alpha, \beta)$ be the multiplier defined by relation (4.13). There exists a complement $W_1$ to $\mathbb{R}\xi$ in $\mathfrak{g}^{G^\xi}$ and two mappings $\rho : \mathbb{R} \times \mathfrak{g}^{G^\xi} \times \mathbb{R}\xi \to \mathfrak{g}^{G^\xi}$ and $\lambda : \mathbb{R} \times W_1 \times \mathbb{R}\xi \to \mathbb{R}$ defined on neighborhoods of the origin such that $\rho(0,0,0) = 0$, $\lambda(0,0) = 0$, and*

$$\frac{\xi + w_0 + \rho(r, \lambda(r, \nu, w_0)\xi + \nu, w_0)}{\Lambda(r, \lambda(r, \nu, w_0)\xi + \nu, w_0 + \rho(r, \lambda(r, \nu, w_0)\xi + \nu, w_0))} = \xi(1 + \lambda(r, \nu, w_0)) + \nu.$$

**Proof.** Let $E : \mathbb{R} \times \mathfrak{g}^{G^\xi} \times \mathfrak{g}^{G^\xi} \to \mathfrak{g}^{G^\xi}$ be the locally defined mapping given by

$$E(r, \alpha, \beta) := \xi + \beta - \Lambda(r, \alpha, \beta)(\xi + \alpha).$$

Note that by Lemma 4.5, $E(0,0,0) = 0$. Now, for each $\beta \in \mathfrak{g}^{G^\xi}$, we have that

$$D_\beta E(0,0,0) \cdot \beta = \left.\frac{d}{dt}\right|_{t=0} (\xi + t\beta - \Lambda(0,0,t\beta)\xi) = \beta - \xi\left(D_\beta \Lambda(0,0,0) \cdot \beta\right).$$

If $\{\xi, \eta_1, \ldots, \eta_p\}$ is a basis of $\mathfrak{g}^{G^\xi}$, then the matrix of the linear map $D_\beta E(0,0,0) : \mathfrak{g}^{G^\xi} \to \mathfrak{g}^{G^\xi}$ in that basis equals:

$$D_\beta E(0,0,0) := \begin{pmatrix} 1 - D_\beta\Lambda(0,0,0) \cdot \xi & -D_\beta\Lambda(0,0,0) \cdot \eta_1 & \cdots & -D_\beta\Lambda(0,0,0) \cdot \eta_p \\ 0 & 1 & \cdots & 0 \\ \vdots & \vdots & \ddots & 0 \\ 0 & 0 & \cdots & 1 \end{pmatrix}.$$

We shall prove that $1 - D_\beta\Lambda(0,0,0) \cdot \xi = 0$. To do this we recall that

$$\Lambda(0,0,\beta) = \frac{\mathbf{d}^2\mathbf{J}^{\xi+\beta}(0)(\psi_2^{-1}(s(0,0,\beta)), \psi_2^{-1}(s(0,0,\beta)))}{\mathbf{d}^2 h(0)(\psi_2^{-1}(s(0,0,\beta)), \psi_2^{-1}(s(0,0,\beta)))}.$$

Therefore,

$$D_\beta\Lambda(0,0,0) \cdot \beta$$
$$= \frac{1}{(\mathbf{d}^2 h(0)(u_0, u_0))^2}\left[\left(\mathbf{d}^2\mathbf{J}^\beta(0)(u_0, u_0) + 2\mathbf{d}^2\mathbf{J}^\xi(0)(D_\beta(\psi_2^{-1} \circ s)(0,0,0) \cdot \beta, u_0))\right)\mathbf{d}^2 h(0)(u_0, u_0)\right.$$
$$\left. - 2\mathbf{d}^2\mathbf{J}^\xi(0)(u_0, u_0)\mathbf{d}^2 h(0)(D_\beta(\psi_2^{-1} \circ s)(0,0,0) \cdot \beta, u_0)\right] = \frac{\mathbf{d}^2\mathbf{J}^\beta(0)(u_0, u_0)}{\mathbf{d}^2 h(0)(u_0, u_0)},$$

where the last equality is a consequence of the fact that $u_0 \in \ker(\mathbf{d}^2(h - \mathbf{J}^\xi)(0))$. Consequently, when we set $\beta = \xi$ in this identity we obtain that $D_\beta\Lambda(0,0,0) \cdot \xi = 1$.

This implies that $W_0 := \ker D_\beta E(0,0,0) = \mathbb{R}\xi$ so by choosing $W_1 := \mathrm{span}\{\eta_1, \ldots, \eta_p\}$ we can write $\mathfrak{g}^{G^\xi} = W_0 \oplus W_1$. Let $\mathbb{P}_{W_0}$ be the projection onto $W_0$. The identity $(\mathbb{I} - \mathbb{P}_{W_0})E(r, \alpha, w_0 + w_1) = 0$ can be solved by the Implicit Function Theorem for $w_1$, which gives us a smooth function $\rho : \mathbb{R} \times \mathfrak{g}^{G^\xi} \times \mathbb{R}\xi \to W_1$ that satisfies

$$(\mathbb{I} - \mathbb{P}_{W_0})E(r, \alpha, w_0 + \rho(r, \alpha, w_0)) \equiv 0. \tag{4.17}$$

Therefore, the solutions of the equation $E(r, \alpha, \beta) = 0$ are in bijective correspondence with the solutions of the scalar equation

$$\mathbb{P}_{W_0}E(r, \alpha, w_0 + \rho(r, \alpha, w_0)) = 0 \tag{4.18}$$

that we will now solve using the Implicit Function Theorem.

***Step 8: Solution of the scalar equation using the Implicit Function Theorem.*** We set

$$g(r, \alpha, w_0) := \mathbb{P}_{W_0} E(r, \alpha, w_0 + \rho(r, \alpha, w_0)) = \xi + w_0 - \Lambda(r, \alpha, w_0 + \rho(r, \alpha, w_0))(\xi + \mathbb{P}_{W_0}\alpha). \quad (4.19)$$

Now, the definition of the function $\rho$ in (4.17) can be rewritten as

$$(\mathbb{I} - \mathbb{P}_{W_0})E(r, \alpha, w_0 + \rho(r, \alpha, w_0)) = \rho(r, \alpha, w_0) - \Lambda(r, \alpha, w_0 + \rho(r, \alpha, w_0))(\mathbb{I} - \mathbb{P}_{W_0})\alpha,$$

which implies that for any value of the parameters $r$ and $w_0$ we have that $\rho(r, 0, w_0) = 0$. Additionally, by implicit differentiation we obtain that $D_\alpha \rho(0,0,0) = \mathbb{I} - \mathbb{P}_{W_0}$. These identities guarantee that $g(0,0,0) = \xi - \Lambda(0,0,0)\xi = 0$ and that

$$\begin{aligned}
D_\alpha g(0,0,0) \cdot \alpha &= -D_\alpha \Lambda(0,0,0) \cdot \alpha - D_\beta \Lambda(0,0,0) \cdot D_\alpha \rho(0,0,0) \cdot \alpha - \mathbb{P}_{W_0} \cdot \alpha \\
&= -D_\alpha \Lambda(0,0,0) \cdot \alpha - \frac{\mathbf{d}^2 \mathbf{J}^{(\mathbb{I}-\mathbb{P}_{W_0})\alpha}(0)(u_0, u_0)}{\mathbf{d}^2 h(0)(u_0, u_0)} - \mathbb{P}_{W_0}\alpha.
\end{aligned} \quad (4.20)$$

We now compute $D_\alpha \Lambda(0,0,0) \cdot \alpha$. Notice that by (4.16) we can write that

$$\Lambda(0, \alpha, 0) = \frac{\mathbf{d}^2 \mathbf{J}^\xi(0)(u(\alpha) + D_{V_0} v_1(0, \alpha) \cdot u(\alpha), u(\alpha) + D_{V_0} v_1(0, \alpha) \cdot u(\alpha))}{\mathbf{d}^2 \mathbf{h}(0)(u(\alpha) + D_{V_0} v_1(0, \alpha) \cdot u(\alpha), u(\alpha) + D_{V_0} v_1(0, \alpha) \cdot u(\alpha))},$$

where $u(\alpha) := \psi_2^{-1}(s(0, \alpha, 0)) \in V_0$. Consequently,

$$\begin{aligned}
D_\alpha \Lambda(0,0,0) \cdot \alpha &= \frac{2\mathbf{d}^2 \mathbf{J}^\xi(0)(D_\alpha u(0) \cdot \alpha + D^2_{V_0, \alpha} v_1(0, 0) \cdot (u_0, \alpha), u_0)\mathbf{d}^2 h(0)(u_0, u_0)}{(\mathbf{d}^2 h(0)(u_0, u_0))^2} \\
&\quad - \frac{2\mathbf{d}^2 h(0)(D_\alpha u(0) \cdot \alpha + D^2_{V_0, \alpha} v_1(0, 0) \cdot (u_0, \alpha), u_0)\mathbf{d}^2 \mathbf{J}^\xi(0)(u_0, u_0)}{(\mathbf{d}^2 h(0)(u_0, u_0))^2}.
\end{aligned}$$

Now, as $u_0 \in \ker(\mathbf{d}^2(h - \mathbf{J}^\xi)(0))$ we have that $\mathbf{d}^2 h(0)(D_\alpha u(0) \cdot \alpha + D^2_{V_0, \alpha} v_1(0, 0) \cdot (u_0, \alpha), u_0) = \mathbf{d}^2 \mathbf{J}^\xi(0)(D_\alpha u(0) \cdot \alpha + D^2_{V_0, \alpha} v_1(0, 0) \cdot (u_0, \alpha), u_0)$ and $\mathbf{d}^2 h(0)(u_0, u_0) = \mathbf{d}^2 \mathbf{J}^\xi(0)(u_0, u_0)$ which substituted in the previous expression implies that $D_\alpha \Lambda(0,0,0) = 0$. Therefore, if in (4.20) we take $\alpha = \xi$ obtain that $D_\alpha g(0,0,0) = -1$ and hence the Implicit Function Theorem guarantees the existence of a function $\lambda: \mathbb{R} \times W_1 \times \mathbb{R}\xi \to \mathbb{R}$ such that $\lambda(0,0,0) = 0$ and

$$g(r, \lambda(r, \nu, w_0)\xi + \nu, w_0) = \mathbb{P}_{W_0} E(r, \lambda(r, \nu, w_0)\xi + \nu, w_0 + \rho(r, \lambda(r, \nu, w_0)\xi + \nu, w_0)) \equiv 0.$$

Finally, the triple $(r, \lambda(r, \nu, w_0)\xi + \nu, w_0 + \rho(r, \lambda(r, \nu, w_0)\xi + \nu, w_0))$ is such that $E(r, \lambda(r, \nu, w_0)\xi + \nu, w_0 + \rho(r, \lambda(r, \nu, w_0)\xi + \nu, w_0)) = 0$ which gives the statement of the lemma for small values of $(r, \nu, w_0)$, since $\Lambda(0,0,0) = 1$. ▼

***Step 9: Closing Arguments.*** By the linearity of the mapping $\mathbf{J}_\alpha^\beta$ in $\beta$, expression (4.13) can be rewritten as

$$\mathbf{d}\mathbf{J}_\alpha^{\frac{\xi + \beta}{\Lambda(r, \alpha, \beta)}}(v_0(r, \alpha, \beta)) = \mathbf{d}H_\alpha(v_0(r, \alpha, \beta)).$$

If we follow the path in the space of parameters $(r, \alpha, \beta)$ given by the functions introduced in Lemma 4.6, that is, $(r, \alpha(r, \nu, w_0), \beta(r, \nu, w_0)) := (r, \lambda(r, \nu, w_0)\xi + \nu, w_0 + \rho(r, \lambda(r, \nu, w_0)\xi + \nu, w_0))$, the above expression becomes

$$\mathbf{d}\mathbf{J}_\alpha^{\xi(1+\lambda(r, \nu, w_0))+\nu}(v_0(r, \alpha(r, \nu, w_0), \beta(r, \nu, w_0))) = \mathbf{d}H_\alpha(v_0(r, \alpha(r, \nu, w_0), \beta(r, \nu, w_0))),$$

or equivalently

$$\nabla_{V_0}(h - \mathbf{J}^{\xi(1+\lambda(r, \nu, w_0))+\nu})(v_0(r, \alpha(r, \nu, w_0), \beta(r, \nu, w_0)) + v_1(v_0(r, \alpha(r, \nu, w_0), \beta(r, \nu, w_0)), \alpha)) = 0.$$

In other words, the pair $(v_0(r, \alpha(r, \nu, w_0), \beta(r, \nu, w_0)), \lambda(r, \nu, w_0)\xi + \nu)$ solves the reduced equation $B(v_0(r, \alpha(r, \nu, w_0), \beta(r, \nu, w_0)), \lambda(r, \nu, w_0)\xi + \nu) = 0$, which implies that the point $v_0(r, \alpha(r, \nu, w_0), \beta(r, \nu, w_0)) +$ $v_1(v_0(r, \alpha(r, \nu, w_0), \beta(r, \nu, w_0)), \lambda(r, \nu, w_0)\xi + \nu) \in V$ is a relative equilibrium of the Hamiltonian vector field $X_h$ with velocity $\xi + \lambda(r, \nu, w_0)\xi + \nu$.

In order to conclude the proof we just need to show that the number of branches predicted in (4.12) coincides with the estimate in the statement of the theorem. Indeed, given that the Lusternik–Schnirelman category takes integer values and the function $H_\alpha$ depends smoothly on $\alpha$, we have that for $\alpha$ small enough

$$\operatorname{Cat}\left(H_\alpha^{-1}(\epsilon)/G^\xi\right) = \operatorname{Cat}\left(H_0^{-1}(\epsilon)/G^\xi\right).$$

The equivariant Morse Lemma, the topologically invariant character of the Lusternik–Schnirelman category, and (4.9) give us that

$$\operatorname{Cat}\left(H_0^{-1}(\epsilon)/G^\xi\right) = \operatorname{Cat}\left(Q^{-1}(\epsilon)/G^\xi\right) = \operatorname{Cat}\left(h|_{V_0}^{-1}(\epsilon)/G^\xi\right), \tag{4.21}$$

where $Q = \mathbf{d}^2 h|_{V_0}(0)$. The estimate (4.3) is a corollary of Proposition 2.2. ∎

## 5 Examples

In this section we illustrate the implementation of Theorem 4.1 with elementary examples that make explicit the procedure suggested by the statement of that result for the study of relative equilibria around symmetric equilibria.

### 5.1 Nonlinearly perturbed spherical pendulum

As it is well known, the spherical pendulum consists of a particle of mass $m$, moving under the action of a constant gravitational field of acceleration $g$, on the surface of a sphere of radius $l$. This system exhibits a circular symmetry obtained when it is rotated around the axis of gravity. The straight down position of the pendulum is a stable equilibrium of the system that is surrounded on each neighboring energy level set by a relative equilibrium. In this example we will use the theorem in the previous section to predict these relative equilibria as well as to show that **they arise in the presence of any $S^1$–invariant nonlinear Hamiltonian perturbation** of the system.

If we use as local coordinates of the configuration space around the downright position the Cartesian coordinates $(x, y)$ of the orthogonal the projection of the sphere on the equatorial plane, the (local) Hamiltonian of this system is:

$$h(x, y, p_x, p_y) = \frac{p_x^2}{2m} + \frac{p_y^2}{2m} - \frac{(xp_x + yp_y)^2}{2ml^2} - mg\sqrt{l^2 - x^2 - y^2} + \varphi(x^2 + y^2, p_x^2 + p_y^2, xp_x + yp_y),$$

where the function $\varphi$ is of order two or higher in all of its variables and encodes the nonlinear perturbation. This system is invariant with respect to the globally Hamiltonian $S^1$–action given by the expression $\Phi_\theta(x, y, p_x, p_y) = (R_\theta(x, y), R_\theta(p_x, p_y))$, where $R_\theta$ denotes a rotation of angle $\theta$. The momentum map $\mathbf{J} : \mathbb{R}^4 \to \mathbb{R}$ associated to this action is given by $\mathbf{J}(x, y, p_x, p_y) = xp_y - yp_x$. The point $(x, y, p_x, p_y) = (0, 0, 0, 0)$ is an equilibrium of the Hamiltonian vector field $X_h$ to which we will apply Theorem 4.1.

Firstly, if $\xi \in \mathbb{R}$ is arbitrary, then

$$\mathbf{d}^2(h - \mathbf{J}^\xi)(0) = \begin{pmatrix} \frac{gm}{l} & 0 & 0 & -\xi \\ 0 & \frac{gm}{l} & \xi & 0 \\ 0 & \xi & \frac{1}{m} & 0 \\ -\xi & 0 & 0 & \frac{1}{m} \end{pmatrix}.$$

Secondly, it is easy to see that $\det(\mathbf{d}^2(h - \mathbf{J}^\xi)(0)) = 0$ iff $\xi = \pm\sqrt{g/l}$. In what follows we will show that on any energy level surrounding the equilibrium there are always two relative equilibria whose velocities are approximately $\pm\sqrt{g/l}$. We will carry out the computations for $\omega := \sqrt{g/l}$. The negative case is completely analogous. It can be verified that

$$V_0^\omega = \operatorname{span}\{(1, 0, 0, m\omega), (0, 1, -m\omega, 0)\}$$

which has a $S^1$–invariant complement $V_1^\omega$ given by

$$V_1^\omega = \operatorname{span}\{(0,0,1,0),(0,0,0,1)\}.$$

We now verify the hypotheses of Theorem 4.1 by writing matricial expression of $\mathbf{d}^2 h(0)|_{V_0}$ using the bases of $V_0^\omega$ and $V_1^\omega$ just given. Indeed,

$$\mathbf{d}^2 h(0)|_{V_0} = \begin{pmatrix} \frac{gm}{l} + m\omega^2 & 0 \\ 0 & \frac{gm}{l} + m\omega^2 \end{pmatrix}$$

which is a positive definite matrix. Let $Q$ be the associated quadratic form. Now, since

$$\Phi_\theta|_{V_0^\omega} = \begin{pmatrix} \cos\theta & -\sin\theta \\ \sin\theta & \cos\theta \end{pmatrix},$$

the $S^1$–action on the circles $Q^{-1}(\epsilon)$ is transitive which forces the functions $j$ defined on it to be necessarily Morse–Bott. Therefore, Theorem 4.1 implies the existence of the relative equilibria that we were looking for.

## 5.2 Coupled oscillators subjected to a magnetic field

The following example provides a situation with higher symmetry than the previous one. We consider the system formed by two identical particles with unit charge and mass $m$ in the $XY$-plane, subjected to identical attractive harmonic forces, to a homogeneous magnetic field perpendicular in direction to the plane of motion $XY$, and to an interaction potential that that will preserve a certain group of symmetries. We will denote by $(q_1, q_2)$ the coordinates of the configuration space of the first particle and by $(q_3, q_4)$ those of the second one. If the magnetic field is induced by the vector potential

$$\mathbf{A}(x,y,z) = \gamma(-y, x, 0),$$

the Hamiltonian function associated to this system is

$$h(\mathbf{q},\mathbf{p}) = \frac{1}{2m}(p_1^2 + p_2^2 + p_3^2 + p_4^2) + \left(\frac{\gamma^2}{2m} + \frac{k}{2}\right)(q_1^2 + q_2^2 + q_3^2 + q_4^2)$$
$$+ \frac{\gamma}{m}(p_1 q_2 - p_2 q_1) + \frac{\gamma}{m}(p_3 q_4 - p_4 q_3) + f(\pi_1, \pi_2, \pi_3, \pi_4), \quad (5.1)$$

where $k$ is a positive constant,

$$\pi_1 = q_1^2 + q_2^2 + q_3^2 + q_4^2$$
$$\pi_2 = p_1^2 + p_2^2 + p_3^2 + p_4^2$$
$$\pi_3 = p_1 q_1 + p_2 q_2 + p_3 q_3 + p_4 q_4$$
$$\pi_4 = p_1 q_2 - p_2 q_1 + p_3 q_4 - p_4 q_3,$$

and $f$ is a function whose order is higher or equal than two in all of its variables. The term involving the function $f$ expresses a non linear interaction between the two particles.

We now study the symmetries of the system. Note that after the assumptions on the interaction function $f$, the system is invariant under the canonical toral action given by the lifted action to the phase space of $R: \mathbb{T}^2 \times \mathbb{R}^4 \to \mathbb{R}^4$, where

$$R((\phi,\psi),\mathbf{q}) = \begin{pmatrix} \cos(\phi)\cos(\psi) & -\cos(\psi)\sin(\phi) & -\cos(\phi)\sin(\psi) & \sin(\phi)\sin(\psi) \\ \cos(\psi)\sin(\phi) & \cos(\phi)\cos(\psi) & -\sin(\phi)\sin(\psi) & -\cos(\phi)\sin(\psi) \\ \cos(\phi)\sin(\psi) & -\sin(\phi)\sin(\psi) & \cos(\phi)\cos(\psi) & -\cos(\psi)\sin(\phi) \\ \sin(\phi)\sin(\psi) & \cos(\phi)\sin(\psi) & \cos(\psi)\sin(\phi) & \cos(\phi)\cos(\psi) \end{pmatrix} \mathbf{q}.$$

and $\mathbf{q} = (q_1, q_3, q_2, q_4)$. The system is also invariant under the transformation

$$\tau \cdot \begin{pmatrix} q_1 \\ q_2 \\ q_3 \\ q_4 \end{pmatrix} = \begin{pmatrix} q_1 \\ q_2 \\ -q_3 \\ -q_4 \end{pmatrix}.$$

The commutation properties of $R$ with the transformation given by $\tau$ make our system $O(2) \times S^1$-invariant. The momentum map $\mathbf{J}: \mathbb{R}^8 \to \mathbb{R}^2$ associated to the toral action is given by the expression

$$\mathbf{J}(\mathbf{q}, \mathbf{p}) = (p_2 q_1 - q_2 p_1 - p_3 q_4 + p_4 q_3, p_3 q_1 - q_3 p_1 - p_2 q_4 + p_4 q_2). \tag{5.2}$$

This system has, for all values of the parameters $\gamma$ and $k$, an equilibrium at the point $(q_1, q_2, q_3, q_4, p_1, p_2, p_3, p_4) = (\mathbf{0}, \mathbf{0})$. We shall use the method described in Theorem 4.1 in order to find the bifurcating relative equilibria from this equilibrium.

Firstly, we find in our particular situation the roots $(\xi_1, \xi_2) \in \mathbb{R}^2$ of equation (4.1), that is,

$$0 = \det\left(\mathbf{d}^2 \left(h - \mathbf{J}^{(\xi_1, \xi_2)}\right)(\mathbf{0}, \mathbf{0})\right)$$
$$= \frac{1}{m^4}\left[k^2 + (\xi_1^2 - \xi_2^2)(4\gamma^2 + 4m\gamma\xi_1 + m^2(\xi_1^2 - \xi_2^2)) - 2k(2\gamma\xi_1 + m(\xi_1^2 + \xi_2^2))\right]^2,$$

which is equivalent to,

$$m^2 \xi_2^4 - 2[(km + \gamma^2) + (m\xi_1 + \gamma)^2]\xi_2^2 + (m\xi_1 + 2\gamma\xi_1 - k)^2 = 0.$$

An analysis of this expression shows that the roots of this equation are given by the pairs $(\xi_1, \xi_2)$ that satisfy any of the four following equalities:

$$\xi_2 = \pm \frac{1}{m}|\xi_1 m + \gamma| \pm \frac{\sqrt{\gamma^2 + km}}{m}. \tag{5.3}$$

We now compute the reduced spaces (the spaces $V_0$ in the notation of Theorem 4.1) associated to the velocities that satisfy (5.3). A detailed study shows that these reduced spaces can be either four or two dimensional. The four dimensional cases correspond to the velocities $\{r_1^+, r_1^-, r_2^+, r_2^-\}$ with corresponding reduced subspaces $\{V_0^{1+}, V_0^{1-}, V_0^{2+}, V_0^{2-}\}$ given by

$$r_1^\pm = \left(\frac{-\gamma \pm \sqrt{km + \gamma^2}}{m}, 0\right), \tag{5.4}$$

$$r_2^\pm = \left(\frac{-\gamma}{m}, \pm\frac{\sqrt{km + \gamma^2}}{m}\right), \tag{5.5}$$

$$V_0^{1\pm} = \mathrm{span}\left\{\left(0, 0, \frac{\pm 1}{\sqrt{km + \gamma^2}}, 0, 0, 0, 0, 1\right), \left(0, 0, 0, \frac{\mp 1}{\sqrt{km + \gamma^2}}, 0, 0, 1, 0\right),\right.$$
$$\left.\left(\frac{\pm 1}{\sqrt{km + \gamma^2}}, 0, 0, 0, 0, 1, 0, 0\right), \left(0, \frac{\mp 1}{\sqrt{km + \gamma^2}}, 0, 0, 1, 0, 0, 0\right)\right\} \tag{5.6}$$

$$V_0^{2\pm} = \mathrm{span}\left\{\left(0, \frac{\pm 1}{\sqrt{km + \gamma^2}}, 0, 0, 0, 0, 0, 1\right), \left(\frac{\pm 1}{\sqrt{km + \gamma^2}}, 0, 0, 0, 0, 0, 1, 0\right),\right.$$
$$\left.\left(0, 0, 0, \frac{\mp 1}{\sqrt{km + \gamma^2}}, 0, 1, 0, 0\right), \left(0, 0, \frac{\mp 1}{\sqrt{km + \gamma^2}}, 0, 1, 0, 0, 0\right)\right\}. \tag{5.7}$$

The two dimensional subspaces correspond to the four one dimensional parameter families of velocities given by

$$r_3^{\pm}(\xi_1) = \left(\xi_1, \pm\left(\frac{1}{m}|\xi_1 m + \gamma| + \frac{\sqrt{\gamma^2 + km}}{m}\right)\right), \quad \xi_1 \in \mathbb{R} \setminus \left\{\frac{-\gamma \pm \sqrt{km + \gamma^2}}{m}, \frac{-\gamma}{m}\right\}, \quad (5.8)$$

$$r_4^{\pm}(\xi_1) = \left(\xi_1, \pm\left(\frac{1}{m}|\xi_1 m + \gamma| - \frac{\sqrt{\gamma^2 + km}}{m}\right)\right), \quad \xi_1 \in \mathbb{R} \setminus \left\{\frac{-\gamma \pm \sqrt{km + \gamma^2}}{m}, \frac{-\gamma}{m}\right\}. \quad (5.9)$$

The associated reduced spaces, that surprisingly do not depend on the parameter $\xi_1$, are given by:

$$V_0^{3\pm} = \text{span}\left\{\left(0, \frac{\pm 1}{\sqrt{km + \gamma^2}}, \frac{-1}{\sqrt{km + \gamma^2}}, 0, \pm 1, 0, 0, 1\right), \left(\frac{\pm 1}{\sqrt{km + \gamma^2}}, 0, 0, \frac{1}{\sqrt{km + \gamma^2}}, 0, \mp 1, 1, 0\right)\right\}, \quad (5.10)$$

$$V_0^{4\pm} = \text{span}\left\{\left(0, \frac{\mp 1}{\sqrt{km + \gamma^2}}, \frac{1}{\sqrt{km + \gamma^2}}, 0, \pm 1, 0, 0, 1\right), \left(\frac{\mp 1}{\sqrt{km + \gamma^2}}, 0, 0, \frac{-1}{\sqrt{km + \gamma^2}}, 0, \mp 1, 1, 0\right)\right\}. \quad (5.11)$$

The quadratic forms $Q_i$ defined as the restrictions $Q_i := \mathbf{d}^2 h(\mathbf{0}, \mathbf{0})|_{V_0^i}$ are given by the expressions:

$$Q_{1\pm} = \frac{2(km + \gamma(\gamma \mp \sqrt{km + \gamma^2}))}{m(km + \gamma^2)} \mathbb{I}_4,$$

$$Q_{2\pm} = \frac{2}{m}\begin{pmatrix} 1 & 0 & 0 & \frac{\gamma}{\sqrt{km+\gamma^2}} \\ 0 & 1 & \frac{-\gamma}{\sqrt{km+\gamma^2}} & 0 \\ 0 & \frac{-\gamma}{\sqrt{km+\gamma^2}} & 1 & 0 \\ \frac{\gamma}{\sqrt{km+\gamma^2}} & 0 & 0 & 1 \end{pmatrix},$$

$$Q_{3\pm} = \frac{4(km + \gamma(\gamma + \sqrt{km + \gamma^2}))}{m(km + \gamma^2)} \mathbb{I}_2,$$

$$Q_{4\pm} = \frac{4(km + \gamma(\gamma - \sqrt{km + \gamma^2}))}{m(km + \gamma^2)} \mathbb{I}_2.$$

The forms $Q_{1\pm}$, $Q_{3\pm}$, and $Q_{4\pm}$ are clearly definite and $Q_{2\pm}$ has as eigenvalues the quantities

$$\frac{2(km + \gamma(\gamma \pm \sqrt{km + \gamma^2}))}{m(km + \gamma^2)}$$

which are always non zero. Hence, $Q_{2\pm}$ is also definite.

The restriction $R|_{V_0^i}$ of the toral action to the reduced spaces $\{V_0^{1+}, V_0^{1-}, V_0^{2+}, V_0^{2-}\}$ always has the same matricial expression if we use as bases the vectors introduced in (5.6) and (5.7), namely,

$$\left(1 + \frac{1}{km + \gamma^2}\right)\begin{pmatrix} \cos(\phi)\cos(\psi) & \cos(\psi)\sin(\phi) & \cos(\phi)\sin(\psi) & \sin(\phi)\sin(\psi) \\ -\cos(\psi)\sin(\phi) & \cos(\phi)\cos(\psi) & -\sin(\phi)\sin(\psi) & \cos(\phi)\sin(\psi) \\ -\cos(\phi)\sin(\psi) & -\sin(\phi)\sin(\psi) & \cos(\phi)\cos(\psi) & \cos(\psi)\sin(\phi) \\ \sin(\phi)\sin(\psi) & -\cos(\phi)\sin(\psi) & -\cos(\psi)\sin(\phi) & \cos(\phi)\cos(\psi) \end{pmatrix}.$$

It can be checked that the eigenvalues of this matrix are given by

$$\left(\frac{1 + km + \gamma^2}{km + \gamma^2}\right)(\cos(\phi \pm \psi) + i\sin(\phi \pm \psi)) \quad \text{and} \quad \left(\frac{1 + km + \gamma^2}{km + \gamma^2}\right)(\cos(\phi \pm \psi) - i\sin(\phi \pm \psi)),$$

which proves that

$$(V_0^{1\pm})^{\mathbb{T}^2} = (V_0^{2\pm})^{\mathbb{T}^2} = \{0\}. \quad (5.12)$$

Additionally,

$$R_{(\phi,\psi)}|_{V_0^{3\pm}} = R_{(\phi,\psi)}|_{V_0^{4\pm}} = \frac{2(1+km+\gamma^2)}{km+\gamma^2} \begin{pmatrix} \cos(\phi \mp \psi) & \sin(\phi \mp \psi) \\ -\sin(\phi \mp \psi) & \cos(\phi \mp \psi) \end{pmatrix},$$

which shows that

$$(V_0^{3\pm})^{\mathbb{T}^2} = (V_0^{4\pm})^{\mathbb{T}^2} = \{0\}, \tag{5.13}$$

that is, the restriction of the toral action to the reduced spaces $\{V_0^{1\pm}, V_0^{2\pm},, V_0^{3\pm}, V_0^{4\pm}\}$ *has trivial fixed point subspaces.*

Finally, it can be verified in a straightforward manner that the restrictions of the quadratic forms $\mathbf{d}^2 \mathbf{J}^{r_i^\pm}(0)$ to the spheres $Q_{i\pm}^{-1}(\epsilon)$ are Morse–Bott functions with respect to the $S^1 \times S^1$–action. Consequently, expressions (5.12) and (5.13) imply that we can use the estimate (4.3) in Theorem 4.1 to conclude that **for each energy level of the system neighboring the origin $(0,0)$ there exist:**

(i) *Eight distinct relative equilibria with respect to the $O(2) \times S^1$ symmetry of the problem that are grouped in four couples; the velocities of the relative equilibria in each couple approach those given by the roots $\{r_1^+, r_1^-, r_2^+, r_2^-\}$ as the energy tends to zero.*

(ii) *Four distinct one parameter families of relative equilibria whose velocities approach those given by the roots $\{r_3^+(\xi_1), r_3^-(\xi_1), r_4^+(\xi_1), r_4^-(\xi_1)\}$ as the energy tends to zero. The parameter $\xi_1$ runs over $\mathbb{R} \setminus \left\{\frac{-\gamma \pm \sqrt{km+\gamma^2}}{m}, \frac{-\gamma}{m}\right\}$.*

## 6 The MGS normal form and the reconstruction equations

In Section 7 we will use the preceding theorems to study the existence of relative equilibria for a Hamiltonian symmetric system in the neighboring energy levels of a stable relative equilibrium that *is not* an equilibrium. The treatment of this problem requires some knowledge of the local geometry and dynamics in symmetric symplectic manifolds, that we will briefly review in this section.

Since this topic has been already introduced already in many other papers we will just briefly sketch the results that we will need in our exposition, and will leave the reader interested in the details consult the original papers [Mar85, GS84]. Regarding the reconstruction equations the reader is encouraged to check with [O98, RWL99, OR02].

Throughout this section we will work with a $G$–Hamiltonian system $(M, \omega, h, G, \mathbf{J})$, where the Lie group $G$ acts in a proper and globally Hamiltonian fashion on the manifold $M$. Let $m$ be a point in $M$ such that $\mathbf{J}(m) = \mu \in \mathfrak{g}^*$ and $G_m$ denotes the isotropy subgroup of the point $m$. We denote by $\mathfrak{g}_\mu$ the Lie algebra of the stabilizer $G_\mu$ of $\mu \in \mathfrak{g}^*$ under the coadjoint action of $G$ on $\mathfrak{g}^*$. We now choose in $\ker T_m \mathbf{J}$ a $G_m$–invariant inner product, $\langle \cdot, \cdot \rangle$, always available by the compactness of $G_m$. Using this inner product we define the ***symplectic normal space*** $V_m$ at $m \in M$ with respect to the inner product $\langle \cdot, \cdot \rangle$, as the orthogonal complement of $T_m(G_\mu \cdot m)$ in $\ker T_m \mathbf{J}$, that is, $\ker T_m \mathbf{J} = T_m(G_\mu \cdot m) \oplus V_m$, where the symbol $\oplus$ denotes orthogonal direct sum. It is easy to verify that $(V_m, \omega(m)|_{V_m})$ is a $G_m$–invariant symplectic vector space.

Recall that by the equivariance of $\mathbf{J}$, the isotropy subgroup $G_m$ of $m$ is a subgroup of $G_\mu$ and therefore $\mathfrak{g}_m = \mathrm{Lie}(G_m) \subset \mathfrak{g}_\mu$. Using again the compactness of $G_m$, we construct an inner product $\langle \cdot, \cdot \rangle$ on $\mathfrak{g}$, invariant under the restriction to $G_m$ of the adjoint action of $G$ on $\mathfrak{g}$. Relative to this inner product we can write the following orthogonal direct sum decompositions $\mathfrak{g} = \mathfrak{g}_\mu \oplus \mathfrak{q}$, and $\mathfrak{g}_\mu = \mathfrak{g}_m \oplus \mathfrak{m}$, for some subspaces $\mathfrak{q} \subset \mathfrak{g}$ and $\mathfrak{m} \subset \mathfrak{g}_\mu$. The inner product also allows us to identify all these Lie algebras with their duals. In particular, we have the dual orthogonal direct sums $\mathfrak{g}^* = \mathfrak{g}_\mu^* \oplus \mathfrak{q}^*$ and $\mathfrak{g}_\mu^* = \mathfrak{g}_m^* \oplus \mathfrak{m}^*$ which allow us to consider $\mathfrak{g}_\mu^*$ as a subspace of $\mathfrak{g}^*$ and, similarly, $\mathfrak{g}_m^*$ and $\mathfrak{m}^*$ as subspaces of $\mathfrak{g}_\mu^*$.

The $G_m$–invariance of the inner product utilized to construct the splittings $\mathfrak{g}_\mu = \mathfrak{g}_m \oplus \mathfrak{m}$ and $\mathfrak{g}_\mu^* = \mathfrak{g}_m^* \oplus \mathfrak{m}^*$, implies that both $\mathfrak{m}$ and $\mathfrak{m}^*$ are $G_m$–spaces using the restriction to them of the $G_m$–adjoint and coadjoint actions, respectively.

The importance of all these objects is in the fact that there is a positive number $r > 0$ such that, denoting by $\mathfrak{m}_r^*$ the open ball of radius $r$ relative to the $G_m$–invariant inner product on $\mathfrak{m}^*$, the manifold $Y_r := G \times_{G_m} (\mathfrak{m}_r^* \times V_m)$ can be endowed with a symplectic structure $\omega_{Y_r}$ with respect to which the left $G$–action $g \cdot [h, \eta, v] = [gh, \eta, v]$ on $Y_r$ is globally Hamiltonian with $\mathrm{Ad}^*$–equivariant momentum map $\mathbf{J}_{Y_r} : Y_r \to \mathfrak{g}^*$ given by $\mathbf{J}_{Y_r}([g, \rho, v]) = \mathrm{Ad}^*_{g^{-1}} \cdot (\mu + \rho + \mathbf{J}_{V_m}(v))$. Moreover, there exist $G$–invariant neighborhoods $U$ of $m$ in $M$, $U'$ of $[e, 0, 0]$ in $Y_r$, and an equivariant symplectomorphism $\phi : U \to U'$ satisfying $\phi(m) = [e, 0, 0]$ and $\mathbf{J}_{Y_r} \circ \phi = \mathbf{J}$. On other words, the twisted product $Y_r$ can be used as a coordinate system in a tubular neighborhood of the orbit $G \cdot m$. This semi–global coordinates are referred to as the **MGS** (Marle–Guillemin–Sternberg) **normal form**.

In what follows we will use the MGS coordinates to compute the equations that describe the dynamics induced by the Hamiltonian vector field corresponding to a $G$–invariant Hamiltonian. These are called the **reconstruction** [O98] or the **bundle** [RWL99] **equations**. Let $h \in C^\infty(Y)^G$ be a $G$–invariant Hamiltonian on $Y$. Our aim is to compute the differential equations that determine the $G$–equivariant Hamiltonian vector field $X_h \in \mathfrak{X}(Y)$ associated to $h$ and characterized by $\mathbf{i}_{X_h} \omega_Y = \mathbf{d}h$.

Since the projection $\pi : G \times \mathfrak{m}^* \times V_m \to G \times_H (\mathfrak{m}^* \times V_m)$ is a surjective submersion, there are always local sections available that we can use to locally express $X_h = T\pi(X_G, X_{\mathfrak{m}^*}, X_{V_m})$, with $X_G$, $X_{\mathfrak{m}^*}$ and $X_{V_m}$ locally defined smooth maps on $Y$ and having values in $TG$, $T\mathfrak{m}^*$ and $TV_m$ respectively. Thus, for any $[g, \rho, v] \in Y$, one has $X_G([g, \rho, v]) \in T_g G$, $X_{\mathfrak{m}^*}([g, \rho, v]) \in T_\rho \mathfrak{m}^* = \mathfrak{m}^*$, and $X_{V_m}([g, \rho, v]) \in T_v V_m = V_m$. Moreover, using the $\mathrm{Ad}_{G_m}$–invariant decomposition of the Lie algebra $\mathfrak{g}$: $\mathfrak{g} = \mathfrak{g}_m \oplus \mathfrak{m} \oplus \mathfrak{q}$, the mapping $X_G$ can be written, for any $[g, \rho, v] \in Y$, as $X_G([g, \rho, v]) = T_e L_g (X_{\mathfrak{g}_m}([g, \rho, v]) + X_{\mathfrak{m}}([g, \rho, v]) + X_{\mathfrak{q}}([g, \rho, v]))$, with $X_{\mathfrak{g}_m}$, $X_{\mathfrak{m}}$, and $X_{\mathfrak{q}}$, locally defined smooth maps on $Y$ with values in $\mathfrak{g}_m$, $\mathfrak{m}$, and $\mathfrak{q}$ respectively. Also, note that since $h \in C^\infty(G \times_H (\mathfrak{m}^* \times V_m))^G$ is $G$–invariant, the mapping $h \circ \pi \in C^\infty(G \times \mathfrak{m}^* \times V_m)^H$ can be understood as a $H$–invariant function that depends only on the $\mathfrak{m}^*$ and $V_m$ variables, that is, $h \circ \pi \in C^\infty(\mathfrak{m}^* \times V_m)^H$.

Using these ideas and the explicit expression of the symplectic form $\omega_{Y_r}$ we can explicitly write down the differential equations that determine the components of $X_h$:

$$X_{\mathfrak{g}_m} = 0 \tag{6.1}$$
$$X_{\mathfrak{q}} = \psi(\rho, v) \tag{6.2}$$
$$X_{\mathfrak{m}} = D_{\mathfrak{m}^*}(h \circ \pi) \tag{6.3}$$
$$X_{V_m} = B_{V_m}^{\sharp}(D_{V_m}(h \circ \pi)) \tag{6.4}$$
$$X_{\mathfrak{m}^*} = \mathbb{P}_{\mathfrak{m}^*}\left(\mathrm{ad}^*_{D_{\mathfrak{m}^*}(h \circ \pi)}\rho + \mathrm{ad}^*_{D_{\mathfrak{m}^*}(h \circ \pi)}\mathbf{J}_{V_m}(v) + \mathrm{ad}^*_{\psi(\rho, v)}(\rho + \mathbf{J}_{V_m}(v))\right). \tag{6.5}$$

**Remark 6.1** The previous equations admit severe simplifications in the presence of various Lie algebraic hypotheses. See [RWL99] for an extensive study. For future reference we will note two particularly important cases:

- The Lie algebra $\mathfrak{g}$ is Abelian: in that case $X_{\mathfrak{m}^*} = X_{\mathfrak{q}} = 0$ at any point.

- The point $\mu \in \mathfrak{g}^*$ is **split** [GLS96], that is, the Lie algebra $\mathfrak{g}_\mu$ of the coadjoint isotropy of $\mu$ admits a $\mathrm{Ad}_{G_\mu}$–invariant complement in $\mathfrak{g}$: in that case the mappings $\eta$ and $\psi$ are identically zero. ◆

**Remark 6.2** The MGS normal form and the reconstruction equations justify why the decision in theorems 3.1 and 4.1 to work with symplectic vector spaces did not imply any loss of generality. Indeed, if in a $G$–Hamiltonian manifold we have an equilibrium $m \in M$ whose isotropy subgroup is $G$, we can locally describe this space around $m$ as $G \times_G V_m \simeq V_m$. In such a situation, the reconstruction equations imply that knowing the dynamics on the $G$–symplectic vector space $V_m$, governed by Hamilton's equations (6.4), is enough to know the dynamics on $G \times_G V_m$ and, therefore, the dynamics in a $G$–invariant neighborhood around $m \in M$. Since theorems 3.1 and 4.1 are local, the claim follows. ◆

# 7 Relative equilibria around a stable relative equilibrium

Our aim in this section is to generalize to relative equilibria, with the help of the MGS normal form and the reconstruction equations, the results that in sections 3 and 4 were proved for equilibria. We start with the generalization of Theorem 3.1. The setup and the notation that will be used coincides with the one introduced in the previous section.

**Theorem 7.1** *Let $(M, \omega, h, G, \mathbf{J})$ be a Hamiltonian $G$–space. Let $m \in M$ be a relative equilibrium of this system with velocity $\xi \in \mathfrak{g}$, such that $H := G_m$, $\mathbf{J}(m) = \mu \in \mathfrak{g}^*$, and $h(m) = 0$. Let $V_m \subset T_m M$ be any symplectic normal space through the point $m$. Suppose that for $V_m$ (and hence for any other symplectic normal space) the following hypotheses are satisfied:*

**(i)** $\mathbf{d}^2 \left( h - \mathbf{J}^{\mathbb{P}_\mathfrak{m} \xi} \right)(m)|_{V_m}$ *is a definite quadratic form.*

**(ii)** $\mathbf{d}^2 \left( \mathbf{J}^{\mathbb{P}_\mathfrak{h} \xi} \right)(m)|_{V_m}$ *is a non degenerate quadratic form.*

**(iii)** *One of the following hypotheses holds:*

  1. *The Lie algebra $\mathfrak{g}$ is Abelian.*
  2. *The Lie algebra $\mathfrak{g}_\mu$ is Abelian and $\mu$ is split.*
  3. $\mathfrak{h} = \mathfrak{g}_\mu$.

*Then for each $\epsilon \in \mathbb{R}$ small enough there are at least*

$$\operatorname{Cat}\left( Q^{-1}(\epsilon)/H^{\mathbb{P}_\mathfrak{h} \xi} \right), \qquad \text{with} \qquad Q(v) = \mathbf{d}^2 \left( h - \mathbf{J}^{\mathbb{P}_\mathfrak{m} \xi} \right)(m)(v,v),\ v \in V_m \qquad (7.1)$$

$H^{\mathbb{P}_\mathfrak{h} \xi}$*–distinct relative equilibria in $h^{-1}(\epsilon)$. The symbol $H^{\mathbb{P}_\mathfrak{h} \xi}$ denotes the adjoint isotropy of the element $\mathbb{P}_\mathfrak{h} \xi \in \mathfrak{h}$ in $H$, and $\operatorname{Cat}$ the Lusternik–Schnirelman category. The projections $\mathbb{P}_\mathfrak{h}$ and $\mathbb{P}_\mathfrak{m}$ are given by the $\operatorname{Ad}_H$–invariant splitting $\mathfrak{g} = \mathfrak{h} \oplus \mathfrak{m} \oplus \mathfrak{q}$ of the Lie algebra $\mathfrak{g}$.*

*If $H^{\mathbb{P}_\mathfrak{h} \xi}$ contains a maximal torus $T^{\mathbb{P}_\mathfrak{h} \xi}$ such that $(V_m)^{T^{\mathbb{P}_\mathfrak{h} \xi}} = \{0\}$, then the estimate (7.1) can be replaced by*

$$\frac{\dim V_m}{2(1 + \dim H^{\mathbb{P}_\mathfrak{h} \xi} - \dim T^{\mathbb{P}_\mathfrak{h} \xi})}. \qquad (7.2)$$

**Remark 7.2** The word *stable* in the title of this section is justified by the fact that condition **(i)** in the statement of Theorem 7.1, along with the existence of a $G_\mu$–invariant inner product on $\mathfrak{g}^*$, with $\mu = \mathbf{J}(m)$, is a sufficient condition [Pat92, LS98, O98, OR99] for the so called $G_\mu$–*stability* of the relative equilibrium $m \in M$. ♦

**Proof** We first verify that the Hessians in the statement are well defined and that the hypotheses on them do not depend on the choice of symplectic normal space $V_m$. As to the first point, it suffices to show that the functions $h - \mathbf{J}^{\mathbb{P}_\mathfrak{m} \xi}$ and $\mathbf{J}^{\mathbb{P}_\mathfrak{h} \xi}$ have $m$ as a critical point. Firstly, since $\xi$ and $\mathbb{P}_\mathfrak{m} \xi$ differ by an element in the Lie algebra of the isotropy of the point $m$ and $\xi$ is a velocity of the relative equilibrium $m$, it follows that $\mathbb{P}_\mathfrak{m} \xi$ is also a velocity and hence we necessarily have $\mathbf{d}(h - \mathbf{J}^{\mathbb{P}_\mathfrak{m} \xi})(m) = 0$. Secondly, since $\mathbb{P}_\mathfrak{h} \xi \in \mathfrak{h}$, the Hamiltonian vector field $X_{\mathbf{J}^{\mathbb{P}_\mathfrak{h} \xi}}$ associated to $\mathbf{J}^{\mathbb{P}_\mathfrak{h} \xi}$ has an equilibrium at the point $m$: $X_{\mathbf{J}^{\mathbb{P}_\mathfrak{h} \xi}}(m) = (\mathbb{P}_\mathfrak{h} \xi)_M(m) = 0$. Therefore, $\mathbf{d}\mathbf{J}^{\mathbb{P}_\mathfrak{h} \xi}(m) = 0$, as required. Regarding the independence of the hypotheses **(i)** and **(ii)** on the choice of symplectic normal space $V_m$, notice that

$$\mathbf{d}^2 \left( h - \mathbf{J}^{\mathbb{P}_\mathfrak{m} \xi} \right)(m)(v, w) = 0, \qquad \mathbf{d}^2 \left( \mathbf{J}^{\mathbb{P}_\mathfrak{h} \xi} \right)(m)(v, w) = 0 \qquad (7.3)$$

whenever $v \in T_m(G_\mu \cdot m)$ and $w \in \ker T_m \mathbf{J}$. Indeed, if we take $v = \eta_M(m)$ with $\eta \in \mathfrak{g}_\mu$, then

$$\mathbf{d}^2(h - \mathbf{J}^{\mathbb{P}_\mathfrak{m} \xi})(m)(v,\, w) = w[X_{\mathbf{J}^\eta}[h - \mathbf{J}^{\mathbb{P}_\mathfrak{m} \xi}]] = w[\{h, \mathbf{J}^\eta\} - \mathbf{J}^{[\mathbb{P}_\mathfrak{m} \xi, \eta]}] = w\left[ \mathbf{J}^{[\mathbb{P}_\mathfrak{m} \xi, \eta]} \right] = 0,$$

where we used the $G_\mu$–invariance of $h$, and that $w \in \ker T_m \mathbf{J}$. The same procedure can be applied to $\mathbf{d}^2 \left( \mathbf{J}^{\mathbb{P}_\mathfrak{h} \xi} \right)(m)$.

Given the local nature of the statement, we can use the MGS coordinates to carry out the proof of the theorem. For simplicity in the exposition we will identify points and maps in $M$ and their counterparts in the MGS coordinates $Y$. Those coordinates can be chosen so that the point $m$ is represented by $[e, 0, 0] \in G \times_H (\mathfrak{m}^* \times V_m)$ and the subset $\Sigma_m := \{e\} \times_H (\{0\} \times V_m) \subset Y$ is a symplectic slice through $m$.

We now verify that $T_m \Sigma_m$ is a symplectic normal space at $m$, that is, $\ker T_m \mathbf{J} = T_m \Sigma_m \oplus T_m(G_\mu \cdot m)$. Indeed, since the canonical projection $\pi : G \times \mathfrak{m}^* \times V_m \longrightarrow G \times_H (\mathfrak{m}^* \times V_m)$ is a surjective submersion, it follows that any vector $v \in T_m M$ can be written as $v = T_{(e,0,0)} \pi(\sigma, \eta, w)$, with some $\sigma \in \mathfrak{g}$, $\eta \in \mathfrak{m}^*$, and $w \in V_m$. In particular, the vectors in $T_m \Sigma_m$ have the form $T_{(e,0,0)} \pi(0, 0, w)$ with $w \in V_m$, and those in $T_m(G_\mu \cdot m)$ look like $T_{(e,0,0)} \pi(\eta, 0, 0)$, with $\eta \in \mathfrak{g}_\mu$. This immediately implies that $T_m \Sigma_m \cap T_m(G_\mu \cdot m) = \{0\}$. At the same time, the equivariance of $\mathbf{J}$ implies that $T_m(G_\mu \cdot m) \subset \ker T_m \mathbf{J}$ and since for any $T_{(e,0,0)} \pi(0, 0, w) \in T_m \Sigma_m$, we have $T_m \mathbf{J}(T_{(e,0,0)} \pi(0, 0, w)) = T_0 \mathbf{J}_{V_m} \cdot w = 0$, it follows that $T_m \Sigma_m \subset \ker T_m \mathbf{J}$. A dimension count shows then that $\ker T_m \mathbf{J} = T_m \Sigma_m \oplus T_m(G_\mu \cdot m)$, as predicted.

Notice that in MGS coordinates the point $m \equiv [e, 0, 0]$ is a relative equilibrium of the Hamiltonian vector field $X_h$ with velocity $\xi$ when

$$X_h(m) = T_{(e,0,0)} \pi(\xi, 0, 0). \tag{7.4}$$

The associated flow is given by $F_t(m) = [\exp t\xi, 0, 0]$.

We now define the function $h_{V_m} \in C^\infty(V_m)^H$ by $h_{V_m}(v) = (h \circ \pi)(0, v)$, for each $v \in V_m$ (as we already said when we introduced the reconstruction equations, a $G$–invariant Hamiltonian in MGS coordinates can be considered as a $H$–invariant function on $\mathfrak{m}^* \times V_m$). Moreover, notice that by (7.4) and the reconstruction equation (6.4)

$$\mathbf{d} h_{V_m}(0) = D_{V_m}(h \circ \pi)(0, 0) = B^\flat_{V_m}(X_{V_m}(0, 0, 0)) = 0,$$

where $B_{V_m} \in \Lambda^2(V_m \times V_m)$ is the Poisson tensor associated to the symplectic form $\omega_{V_m} := \omega|_{V_m}$ and $B^\flat_{V_m} : V_m \to V_m$ is the associated linear map. Also, for any $v, w \in V_m$:

$$\mathbf{d}^2(h - \mathbf{J}^{\mathbb{P}_\mathfrak{m} \xi})([e, 0, 0])(T_{(e,0,0)} \pi(0, 0, v), T_{(e,0,0)} \pi(0, 0, w)) = \left.\frac{d}{dt}\right|_{t=0} \left.\frac{d}{ds}\right|_{s=0} (h - \mathbf{J}^{\mathbb{P}_\mathfrak{m} \xi})([e, 0, tv + sw])$$

$$= \mathbf{d}^2 h_{V_m}(0)(v, w) - \left.\frac{d}{dt}\right|_{t=0} \langle T_{tv} \mathbf{J}_{V_m} \cdot w, \mathbb{P}_\mathfrak{m} \xi \rangle = \mathbf{d}^2 h_{V_m}(0)(v, w),$$

since $T_{tv} \mathbf{J}_{V_m} \cdot w \in \mathfrak{h}^*$ for any $t$. Therefore, hypothesis **(i)** implies that $\mathbf{d}^2 h_{V_m}(0)$ is a definite quadratic form. Analogously, we can show that for any $v, w \in V_m$:

$$\mathbf{d}^2(\mathbf{J}^{\mathbb{P}_\mathfrak{h} \xi})([e, 0, 0])(T_{(e,0,0)} \pi(0, 0, v), T_{(e,0,0)} \pi(0, 0, w))$$

$$= \left.\frac{d}{dt}\right|_{t=0} \left.\frac{d}{ds}\right|_{s=0} (\mathbf{J}^{\mathbb{P}_\mathfrak{h} \xi})([e, 0, tv + sw]) = \mathbf{d}^2 \mathbf{J}^{\mathbb{P}_\mathfrak{h} \xi}_{V_m}(0)(v, w),$$

which by hypothesis **(ii)** implies that $\mathbf{d}^2 \mathbf{J}^{\mathbb{P}_\mathfrak{h} \xi}_{V_m}(0)$ is a nondegenerate quadratic form.

If we now apply Theorem 3.1 to the equilibrium that the system $(V_m, \omega_{V_m}, h_{V_m}, H, \mathbf{J}_{V_m})$ has at the origin we obtain at least

$$\operatorname{Cat}\left(Q^{-1}(\epsilon)/H^{\mathbb{P}_\mathfrak{h} \xi}\right), \quad \text{with} \quad Q(v) = \mathbf{d}^2 h_{V_m}(0)(v, v) = \mathbf{d}^2\left(h - \mathbf{J}^{\mathbb{P}_\mathfrak{m} \xi}\right)(m)(v, v), \; v \in V_m \tag{7.5}$$

$H$–relative equilibria for that system whose velocities are a real multiple of $\mathbb{P}_\mathfrak{h} \xi$.

In the rest of the proof we will see that the hypotheses in the assumption **(iii)** of the statement allow us to use these $H$–relative equilibria to construct $G$–relative equilibria of the original system. Suppose that we are in the first two cases considered in the hypothesis **(iii)**, that is, either $\mathfrak{g}_\mu$ is Abelian and $\mu$ split or $\mathfrak{g}$ is Abelian. Having in mind what we said in Remark 6.1 and the reconstruction equation (6.5), we realize that $X_{\mathfrak{m}^*} = 0$ at any point and therefore if $v \in V_m$ is one of the $H$–relative equilibria of $(V_m, \omega_{V_m}, h_{V_m}, H, \mathbf{J}_{V_m})$ predicted by (7.5), the point $[e, 0, v]$ is necessarily a $G$–relative equilibrium of the original system.

Finally, if $\mathfrak{h} = \mathfrak{g}_\mu$, then $\mathfrak{m} = 0$ necessarily and therefore all the points of the form $[e, v]$, with $v \in V_m$, one of the $H$–relative equilibria predicted by (7.5), are $G$–relative equilibria of the original system. ∎

We finish with the generalization to relative equilibria of Theorem 4.1.

**Theorem 7.3** *Let $(M, \omega, h, G, \mathbf{J})$ be a Hamiltonian $G$–space. Let $m \in M$ be a relative equilibrium of this system with velocity $\xi \in \mathfrak{g}$, such that $H := G_m$, $\mathbf{J}(m) = \mu \in \mathfrak{g}^*$, and $h(m) = 0$. Let $V_m \subset T_m M$ be any symplectic normal space through the point $m$. Let $\eta \in \mathfrak{h}$ be a root of the polynomial equation:*

$$\det \left( \mathbf{d}^2 \left( h - \mathbf{J}^{\mathbb{P}_\mathfrak{m} \xi + \eta} \right)(m)|_{V_m} \right) = 0.$$

*Define the subspace $V_0 \subset V_m$ by*

$$V_0 := \ker \left( \mathbf{d}^2 \left( h - \mathbf{J}^{\mathbb{P}_\mathfrak{m} \xi + \eta} \right)(m)|_{V_m} \right).$$

*Suppose that the following hypotheses are satisfied:*

**(i)** $\mathbf{d}^2 \left( h - \mathbf{J}^{\mathbb{P}_\mathfrak{m} \xi} \right)(m)|_{V_0}$ *is a definite quadratic form.*

**(ii)** *Let $Q(v) := \mathbf{d}^2 \left( h - \mathbf{J}^{\mathbb{P}_\mathfrak{m} \xi} \right)(m)(v, v)$, $v \in V_0$, $\| \cdot \|$ be the norm on $V_0$ associated to $Q$, $l = \dim V_0$, and $S^{l-1}$ be the unit sphere in $V_0$ defined by the norm $\| \cdot \|$. The function $j^\eta \in C^\infty(S^{l-1})$ defined by $j^\eta(u) := \frac{1}{2} \mathbf{d}^2 \mathbf{J}^\eta(0)(u, u)$, is Morse–Bott with respect to the $H^\eta$–action on $S^{l-1}$.*

**(iii)** *One of the following hypotheses holds:*

   1. *The Lie algebra $\mathfrak{g}$ is Abelian.*
   2. *The Lie algebra $\mathfrak{g}_\mu$ is Abelian and $\mu$ is split.*
   3. $\mathfrak{h} = \mathfrak{g}_\mu$.

*Then for each $\epsilon \in \mathbb{R}$ small enough there are at least*

$$\operatorname{Cat}\left( Q^{-1}(\epsilon)/H^\eta \right) \tag{7.6}$$

*$H^\eta$–distinct relative equilibria in $h^{-1}(\epsilon)$. The symbol $H^\eta$ denotes the adjoint isotropy of the element $\eta \in \mathfrak{h}$ in $H$ and $\operatorname{Cat}$ the Lusternik–Schnirelman category. The projection $\mathbb{P}_\mathfrak{m}$ is given by the $\operatorname{Ad}_H$–invariant splitting $\mathfrak{g} = \mathfrak{h} \oplus \mathfrak{m} \oplus \mathfrak{q}$ of the Lie algebra $\mathfrak{g}$.*

*If the subgroup $H^\eta$ contains a maximal torus $T^\eta$ such that the subspace $(V_0)^{T^\eta}$ of $T^\eta$–fixed point vectors in $V_0$ is trivial, that is $(V_0)^{T^\eta} = \{0\}$, then the estimate (7.6) can be replaced by*

$$\frac{\dim V_0}{2(1 + \dim H^\eta - \dim T^\eta)}. \tag{7.7}$$

**Proof** It suffices to reproduce the *modus operandi* followed in the proof of Theorem 7.1, this time invoking Theorem 4.1 once the $H$–invariant Hamiltonian dynamical system $(V_m, \omega_{V_m}, h_{V_m})$ has been constructed and the hypotheses in the statement have been used. ∎

**Acknowledgments.** We thank R. Cushman and J. Montaldi for several illuminating discussions during the elaboration of this work. This research was partially supported by the European Commission and the Swiss Federal Government through funding for the Research Training Network *Mechanics and Symmetry in Europe* (MASIE). The support of the Swiss National Science Foundation is also acknowledged.